\documentclass[12pt,twoside]{article}
\usepackage{ajc}

\Volume{42}\Year{2008}\firstpage{99}\lastpage{114}\Received{26 July 2007}\Revised{13 Nov 2007}
\runninghead{Perfect domination in regular grid graphs}
\shortauthors{Italo J. Dejter}
\newtheorem{thm}{Theorem} 

\newtheorem{cor}[thm]{Corollary}

\newcommand{\Z}{\hbox{\bf Z}}
\newcommand{\R}{\hbox{\bf R}}

\newcommand{\beeq}{\begin{eqnarray*}}
\newcommand{\eneq}{\end{eqnarray*}}
\newcommand{\proof}{\noindent {\it Proof.\hspace{4mm}}}

\newcommand{\qfd}{\hfill $\fbox{}$\vspace{4mm}}\def\newpic#1{%
\def\emline##1##2##3##4##5##6{%
\put(##1,##2){\special{em:point #1##3}}%
\put(##4,##5){\special{em:point #1##6}}%
\special{em:line #1##3,#1##6}}}
\newpic{}
\def\emline#1#2#3#4#5#6{%
\put(#1,#2){\special{em:moveto}}%
\put(#4,#5){\special{em:lineto}}}
\def\newpic#1{}
\newcommand\ZZ{{{\rm Z}\kern-.28em{\rm Z}}}

\title{Perfect domination in regular grid graphs}
\author{Italo J. Dejter
\\ University of Puerto Rico \\ Rio Piedras, PR 00931-3355 \\ ijdejter@uprrp.edu
}
\date{}

\begin{document}
\maketitle

\begin{abstract}
We show there is an uncountable number of parallel total perfect
codes in the integer lattice graph ${\Lambda}$ of $\R^2$. In
contrast, there is just one 1-perfect code in ${\Lambda}$ and one
total perfect code in ${\Lambda}$ restricting to total perfect codes
of rectangular grid graphs (yielding an asymmetric, Penrose, tiling
of the plane). We characterize all cycle products $C_m\times C_n$
with parallel total perfect codes, and the $d$-perfect and total
perfect code partitions of ${\Lambda}$ and $C_m\times C_n$, the
former having as quotient graph the undirected Cayley graphs of
$\Z_{2d^2+2d+1}$ with generator set $\{1,2d^2\}$. For $r>1$,
generalization  for 1-perfect codes is provided in the integer
lattice of $\R^r$ and in the products of $r$ cycles, with partition
quotient graph $K_{2r+1}$ taken as the undirected Cayley graph of
$\Z_{2r+1}$ with generator set $\{1,\ldots,r\}$.
\end{abstract}

\section{Introduction}

As in \cite{W}, a vertex subset $S$ in a graph $G$ is said to be a
{\it perfect dominating set} (PDS) in $G$ if each vertex of the
complementary graph $G\setminus S$ of $S$ in $G$ is adjacent to just
one element of $S$. In that case, if the induced components of $S$
are 1-cubes, (respectively 0-cubes), then $S$ is said to be a {\it
total perfect code} \cite{KG}, (respectively a {\it 1-perfect code}
\cite{K}, or {\it efficient dominating set} \cite{BBS}).

vspace*{5mm}

The NP-completeness of finding a 1-perfect code in $G$ and that of
finding a minimal PDS in a planar graph were established
respectively in \cite{BBS,K} and in Sections 3 and 4 of \cite{FH},
even if its induced components are $i$-cubes with $i\leq 1$.

The {\it integer lattice graph} ${\Lambda}$ of $\R^2$ is the graph
with vertex set $\{(i,j):i,j\in\Z\}$ and such that any two vertices
in ${\Lambda}$ are adjacent if and only if their Euclidean distance
is 1. ${\Lambda}$ and its subgraphs are represented orthogonally,
their vertical paths from left to right for increasing indices
$i=0,1,2\ldots,m-1$ and its horizontal paths downward for increasing
indices $j=0,1,2\dots,n-1$. If a total perfect code in ${\Lambda}$
has its induced components as pairwise parallel 1-cubes in
${\Lambda}$, then it is said to be {\it parallel}.

We represent a perfect dominating set $S$ in ${\Lambda}$
by setting the vertices of $S$ as black dots at their locations and by tracing only the
edges between vertices in ${\Lambda}\setminus S$ (by means of unit-length solid segments),
thus avoiding (or deleting) those edges in ${\Lambda}$ incident to $S$.
This way, we represent not only $S$ but also ${\Lambda}\setminus S$,
as in Figure 1, that accompanies Theorem 1. Notice that $G\setminus S$ is
3-regular, for any PDS $S$ in a 4-regular graph $G$, like $G={\Lambda}$.

We show in Section 2 that there is no algorithmic characterization of parallel total perfect
codes in ${\Lambda}$.
In fact, there is an uncountable number of such codes, in one-to-one correspondence
with the doubly infinite binary sequences.
This allows to characterize, in Section 3, all cycle products $C_m\times C_n$ in which
parallel total perfect codes exist, and in Sections 4 and 5,
the $d$-perfect and total perfect code partitions of ${\Lambda}$ and $C_m\times C_n$.
We also show that the quotient graphs of these
$d$-perfect code partitions are Cayley graphs of $\Z_{2d^2+2d+1}$ with
generator sets $\{1,2d^2\}$.
In contrast with the mentioned uncountability, seen as one of perfect dominating sets,
not only there is just one 1-perfect code in ${\Lambda}$ but, as a result of a
characterization of grid graphs containing total perfect codes
due to Klostermeyer and Goldwasser, there exists only one total perfect code
in ${\Lambda}$ that restricts to total perfect codes
of rectangular grid graphs, its complementary graph in ${\Lambda}$ yielding
an asymmetric tiling of the plane, like the Penrose tiling, (Section 6).
Section 7 considers some special cases, including a correction to a result of \cite{LS}.
In Section 8, the generalization of  the results for 1-perfect codes is provided
in dimensions $r>1$, with partition quotient graph $K_{2r+1}$ as undirected
Cayley graph of $\Z_{2r+1}$ with generator set $\{1,\ldots,r\}$.

\begin{figure}
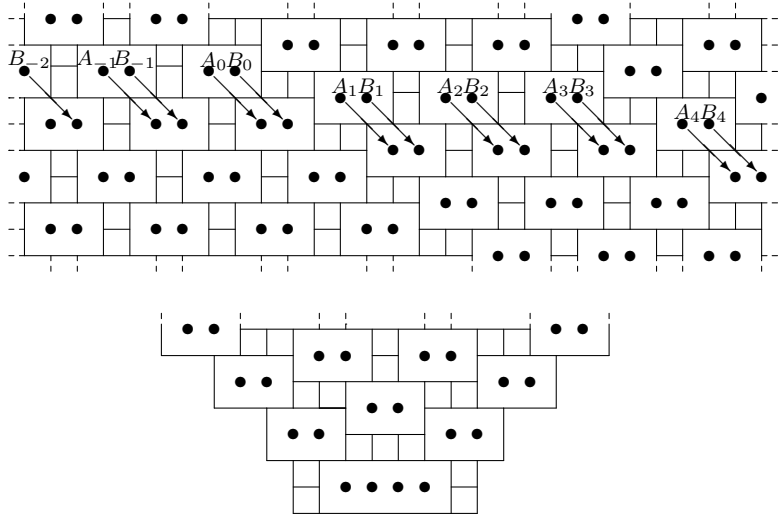

\unitlength=0.70mm \special{em:linewidth 0.4pt}
\linethickness{0.4pt}

\caption{Illustrations for Theorem 1}
\end{figure}

\section{Total perfect codes in ${\Lambda}$}

\begin{thm}
The family of PDSs in ${\Lambda}$ whose induced components are
parallel, horizontal, 1-cubes is in one-to-one correspondence with
the set of doubly infinite $\{0,1\}$-sequences.
\end{thm}

\proof Let $(\ldots,a_{-i},\ldots,a_0, \ldots, a_i,\ldots)$ be a
doubly infinite $\{0,1\}$-sequence.

Define vertices $A_n,B_n$ in ${\Lambda}$ by means of
$$A_{-n}=\Sigma_{i=-1}^{-n}-(4+a_i,a_i),\,\, A_0=(0,0),\,\,
                         A_n=\Sigma_{i=1}^n(4+a_i,a_i),$$
for $n>0$; and $B_n=A_n+(1,0)$, for $n\in\Z$. Then,
$$\cup_{j=-\infty}^{\infty}s^j(\cup_{i=-\infty}^{\infty}\{A_i,B_i\})$$
is a PDS in ${\Lambda}$, where $s^j(x,y)=(x+2j,y+2j)$, with $j\in\Z$.
This yields the claimed correspondence.
For example,
$$(\ldots,a_{-2},a_{-1},a_0,a_1,a_2,a_3,a_4,a_5,\ldots)=(\ldots, 0,0,1,0,0,1,1, \ldots)$$
is sent into a PDS $S$ as partially
represented on the upper part of Figure 1,
where only the edges in ${\Lambda}\setminus S$ are shown and the arrows indicate
assignment of some vertices $A_n$ and $B_n$, with $|n|$ small, via $s=s^1$.

Now, PDSs whose induced components are
parallel, horizontal, 1-cubes
can only be images of $\{0,1\}$-sequences through the just
introduced one-to-one correspondence, for otherwise a path of length 3 would be found
as an induced component of such a PDS,
as in the bottom of Figure 1, a contradiction. In particular, we conclude that there are
no less PDSs in ${\Lambda}$ than doubly infinite binary sequences, whose cardinality
surpasses $\aleph_0$.
\qfd

\begin{cor}
There is no algorithmic characterization for PDSs in ${\Lambda}$.
\end{cor}

\proof
Because of the uncountability of total perfect codes in ${\Lambda}$ found in
Theorem 1, it is clear that there cannot be such algorithmic characterization. \qfd

Given a PDS $S$ in ${\Lambda}$, the graph ${\Lambda}\setminus S$ has chordless
cycles delimiting rectangles of areas at least 4,
that we call {\it rooms}, and maximal connected unions of 4-cycles arranged either
horizontally or vertically into rectangles,
that we call {\it ladders}, a possible case of which is a
single 4-cycle bordered by rooms.
The totality of ordered pairs formed by the
horizontal and vertical dimensions, (widths and heights), of the
rectangles spanned by these rooms and ladders can be presented
in an array of integer pairs that we call the {\it PDS-array} associated to $S$.
For example, the PDS partially depicted in the upper part of Figure 1 has a
PDS-array correspondingly containing the following disposition of ordered pairs
of (one-digit) integers (with deleted parentheses and commas):
$$\begin{array}{cccccccccccccccc}
\ldots & 12 & 32 & 12 & 32 & 21 & 32 & 12 & 32 & 12 & 32 & 21 & 32 & 21 & 32 & \ldots \\
\ldots & 32 & 12 & 32 & 12 & 32 & 21 & 32 & 12 & 32 & 12 & 32 & 21 & 32 & 21 & \ldots \\
\ldots & 12 & 32 & 12 & 32 & 12 & 32 & 21 & 32 & 12 & 32 & 12 & 32 & 21 & 32 & \ldots \\
\ldots & 32 & 12 & 32 & 12 & 32 & 32 & 32 & 21 & 32 & 12 & 32 & 12 & 32 & 21 & \ldots \\
\ldots & 12 & 32 & 12 & 32 & 12 & 12 & 12 & 32 & 21 & 32 & 12 & 32 & 12 & 32 & \ldots \\
\end{array}$$
with the pairs $32=(3,2)$ in the second and third rows representing the rooms
containing the vertices $A_i,B_i$ in the figure and its destination vertices $s(A_i),$
$s(B_i)$, respectively, for $i=-2,\ldots,4$. The pairs in a PDS-array giving the dimensions of
a room, (ladder), are called {\it room pairs}, ({\it ladder pairs}).

Let TPC stand for total perfect code.
Let $\Psi$ denote the correspondence in Theorem 1 that assigns to each doubly infinite
binary sequence $(\ldots,a_i,$ $\ldots,a_0,$ $\ldots,a_i\ldots)$
a {\it parallel} TPC (PTPC) in ${\Lambda}$,
meaning that its induced 1-cubes are pairwise parallel.

\begin{cor}
The family of TPCs $S$ in ${\Lambda}$ having only ladder pairs of the
form $(2,1)=21$ is in one-to-one correspondence with
the set of doubly infinite $\{0,1\}$-sequences.
\end{cor}

\proof
Consider the translated lattice graph
${\Lambda}'={\Lambda}+(\frac{1}{2},\frac{1}{2})$.
Given a doubly infinite binary sequence $A=(\ldots,a_i,\ldots,a_0,\ldots,a_i\ldots)$,
we construct a TPC $\Psi'(A)$
in ${\Lambda}'$ as in the statement of the corollary
by selecting its component vertices
(in ${\Lambda}'$, instead of ${\Lambda}$)
as the barycenters of the unit-area squares delimiting the 4-cycles of
${\Lambda}\setminus \Psi(A)$. Then, by taking
$\Psi''(A)=\Psi'(A)-(\frac{1}{2},\frac{1}{2})$, we get a one-one correspondence as required.
\qfd

The PDS-array partially presented above in relation to the example
of Figure 1, for a PTPC $\Psi(\ldots,0,0,1,0,0,1,1,\ldots)$, has a
counterpart due to the argument of Corollary 3, which is the PDS-array for
the TPC
$$\Psi''(\ldots,0,0,1,0,0,1,1,\ldots),$$ correspondingly representable as follows:

$$\begin{array}{cccccccccccccccc}
\ldots & 23 & 21 & 23 & 21 & 32 & 21 & 23 & 21 & 23 & 21 & 32 & 21 & 32 & 21 & \ldots \\
\ldots & 21 & 23 & 21 & 23 & 21 & 32 & 21 & 23 & 21 & 23 & 21 & 32 & 21 & 32 & \ldots \\
\ldots & 23 & 21 & 23 & 21 & 23 & 21 & 32 & 21 & 23 & 21 & 23 & 21 & 32 & 21 & \ldots \\
\ldots & 21 & 23 & 21 & 23 & 21 & 21 & 21 & 32 & 21 & 23 & 21 & 23 & 21 & 32 & \ldots \\
\ldots & 23 & 21 & 23 & 21 & 23 & 23 & 23 & 21 & 32 & 21 & 23 & 21 & 23 & 21 & \ldots \\
\end{array}$$
obtained from the previously given PDS-array
by subtracting 1 from each one of the two component integers of each room pair and
adding 1 to each one of the two component integers of each ladder pair.

We may say that the relation between $\Psi(A)$ and $\Psi'(A)$ is one of {\it room-ladder
duality}. This applies, more generally,
to all PDSs $S$ in ${\Lambda}$
whose complementary graphs ${\Lambda}\setminus S$ in ${\Lambda}$ have PDS-arrays with
room pairs $(r,s)$ satisfying only $min\{r,s\}=2$. This extends to PDSs in
quotient graphs in ${\Lambda}$ of the form $C_m\times C_n$, where $min\{m,n\}\geq 3$.

Theorem 1 may have been expressed for PDSs whose induced components are parallel, vertical,
1-cubes. Corollary 3 may have been expressed for TPCs having only ladder pairs of the form
$(1,2)=12$.

\section{PTPCs in cycle products}

By restricting $\Psi$ to periodic doubly infinite
binary sequences, we find the following result.
A binary $n$-tuple has {\it weight} $k$ if and only if it has exactly $k$ unit entries.

\begin{thm} A periodic
doubly infinite binary sequence with period $B=(b_1,b_2,\ldots,$ $b_n)$ of
weight $k$ is sent by $\Psi$ into a PTPC $\Psi(B)$ in ${\Lambda}$ that
covers a PTPC $\Phi(B)$ in the Cartesian product $C_p\times C_p$, where
$p=4n$ if $k$ is even, and $p=8n$ if $k$ is odd.
This establishes a one-one correspondence $\Phi$ from the irreducible 
periodic doubly infinite binary
sequences onto the irreducible horizontal PTPCs in Cartesian products of
two cycles.
\end{thm}

\proof The irreducible periods $B=(b_1,b_2,\ldots,b_n)$
of periodic doubly infinite binary sequences can be listed
lexicographically and by their nondecreasing lengths, as follows:
0, 1, 01, 001, 011, 0001, 0011, 0111, 00001, 00011, 00101, 00111, 01011, 01111, $\ldots,$
where periods whose images through $\Psi$ are horizontal
PTPCs $\Psi(B)$ in ${\Lambda}$ already obtained from previously listed
periods are avoided. This yields two subsequences according to
weight parity, namely: 0, 011, 0011, 00011, 00101, 01111, $\ldots$ and
1, 01, 001, 0001, 0111, 00001, 00111, 01011, $\ldots,$ with respective
sequences of values of $p$ as in the statement of the theorem as follows:
$p=$ 4, 12, 16, 20, 20, 20, $\ldots$ and $p=$ 8, 16, 24, 32, 32, 40, 40, 40, $\ldots.$

In the PDS-array associated to one such $\Psi(B)$, there are only two types of
horizontally contiguous pairs of entries (which are integer pairs)
formed by a leftmost room pair and a rightmost ladder pair,
namely: $(32,12)$ and $(32,21)$. For example, Figure 2 shows the first stage of the
procedure in the proof of Theorem 1 applied to the irreducible period
$B=00011101$, producing concretely
the {\it fundamental tile} ${\mathcal T}(B)$ associated to $B$,
before extending it periodically to its (upper-)left and to its (lower-)right sides
and applying subsequently the transformation $s$ of Theorem 1 and its powers
in order to exhibit $\Psi(B)$.
Observe that each null entry in this $B$ has associated a
pair $(32,12)$
and each non-null entry in it has associated a (descending) pair $(32,21)$.

\begin{figure}
\unitlength=0.65mm
\special{em:linewidth 0.4pt}
\linethickness{0.4pt}
\begin{picture}(207.43,26.89)
\emline{26.43}{25.89}{1}{26.43}{15.89}{2}
\emline{26.43}{15.89}{3}{41.43}{15.89}{4}
\emline{41.43}{15.89}{5}{41.43}{25.89}{6}
\emline{41.43}{25.89}{7}{26.43}{25.89}{8}
\emline{26.43}{25.89}{9}{41.43}{25.89}{10}
\emline{41.43}{25.89}{11}{46.43}{25.89}{12}
\emline{46.43}{25.89}{13}{46.43}{15.89}{14}
\emline{46.43}{15.89}{15}{41.43}{15.89}{16}
\emline{41.43}{20.89}{17}{46.43}{20.89}{18}
\emline{41.43}{25.89}{19}{26.43}{25.89}{20}
\put(31.43,20.89){\circle*{2.00}}
\put(36.43,20.89){\circle*{2.00}}
\emline{46.43}{25.89}{21}{46.43}{15.89}{22}
\emline{46.43}{15.89}{23}{61.43}{15.89}{24}
\emline{61.43}{15.89}{25}{61.43}{25.89}{26}
\emline{61.43}{25.89}{27}{46.43}{25.89}{28}
\emline{61.43}{25.89}{29}{66.43}{25.89}{30}
\emline{66.43}{15.89}{31}{61.43}{15.89}{32}
\emline{61.43}{20.89}{33}{66.43}{20.89}{34}
\put(51.43,20.89){\circle*{2.00}}
\put(56.43,20.89){\circle*{2.00}}
\emline{66.43}{25.89}{35}{66.43}{15.89}{36}
\emline{66.43}{15.89}{37}{81.43}{15.89}{38}
\emline{81.43}{15.89}{39}{81.43}{25.89}{40}
\emline{81.43}{25.89}{41}{66.43}{25.89}{42}
\emline{81.43}{25.89}{43}{86.43}{25.89}{44}
\emline{86.43}{25.89}{45}{86.43}{15.89}{46}
\emline{86.43}{15.89}{47}{81.43}{15.89}{48}
\emline{81.43}{20.89}{49}{86.43}{20.89}{50}
\put(71.43,20.89){\circle*{2.00}}
\put(76.43,20.89){\circle*{2.00}}
\emline{86.43}{15.89}{51}{101.43}{15.89}{52}
\emline{101.43}{15.89}{53}{101.43}{25.89}{54}
\emline{101.43}{25.89}{55}{86.43}{25.89}{56}
\emline{106.43}{15.89}{57}{101.43}{15.89}{58}
\emline{101.43}{20.89}{59}{106.43}{20.89}{60}
\put(91.43,20.89){\circle*{2.00}}
\put(96.43,20.89){\circle*{2.00}}
\emline{106.43}{20.89}{61}{106.43}{15.89}{62}
\emline{106.43}{20.89}{63}{111.43}{20.89}{64}
\emline{111.43}{20.89}{65}{111.43}{15.89}{66}
\emline{111.43}{15.89}{67}{106.43}{15.89}{68}
\emline{126.43}{20.89}{69}{126.43}{10.89}{70}
\emline{111.43}{10.89}{71}{111.43}{15.89}{72}
\put(116.43,15.89){\circle*{2.00}}
\put(121.43,15.89){\circle*{2.00}}
\emline{131.43}{10.89}{73}{126.43}{10.89}{74}
\emline{126.43}{15.89}{75}{131.43}{15.89}{76}
\emline{131.43}{15.89}{77}{131.43}{10.89}{78}
\emline{131.43}{15.89}{79}{136.43}{15.89}{80}
\emline{136.43}{15.89}{81}{136.43}{10.89}{82}
\emline{136.43}{10.89}{83}{131.43}{10.89}{84}
\emline{151.43}{15.89}{85}{151.43}{5.89}{86}
\emline{136.43}{5.89}{87}{136.43}{10.89}{88}
\put(141.43,10.89){\circle*{2.00}}
\put(146.43,10.89){\circle*{2.00}}
\emline{156.43}{5.89}{89}{151.43}{5.89}{90}
\emline{151.43}{10.89}{91}{156.43}{10.89}{92}
\emline{156.43}{10.89}{93}{156.43}{5.89}{94}
\emline{156.43}{10.89}{95}{161.43}{10.89}{96}
\emline{161.43}{10.89}{97}{161.43}{5.89}{98}
\emline{161.43}{5.89}{99}{156.43}{5.89}{100}
\emline{161.43}{0.89}{101}{161.43}{5.89}{102}
\put(166.43,5.89){\circle*{2.00}}
\put(171.43,5.89){\circle*{2.00}}
\emline{111.43}{20.89}{103}{126.43}{20.89}{104}
\emline{111.43}{10.89}{105}{126.43}{10.89}{106}
\emline{136.43}{15.89}{107}{151.43}{15.89}{108}
\emline{136.43}{5.89}{109}{151.43}{5.89}{110}
\emline{161.43}{10.89}{111}{176.43}{10.89}{112}
\emline{176.43}{10.89}{113}{176.43}{0.89}{114}
\emline{176.43}{0.89}{115}{161.43}{0.89}{116}
\emline{176.43}{0.89}{117}{181.43}{0.89}{118}
\emline{181.43}{0.89}{119}{181.43}{10.89}{120}
\emline{181.43}{10.89}{121}{176.43}{10.89}{122}
\emline{176.43}{5.89}{123}{181.43}{5.89}{124}
\emline{196.43}{10.89}{125}{196.43}{0.89}{126}
\emline{201.43}{0.89}{127}{196.43}{0.89}{128}
\emline{196.43}{5.89}{129}{201.43}{5.89}{130}
\emline{201.43}{5.89}{131}{201.43}{0.89}{132}
\emline{201.43}{5.89}{133}{206.43}{5.89}{134}
\emline{206.43}{5.89}{135}{206.43}{0.89}{136}
\emline{206.43}{0.89}{137}{201.43}{0.89}{138}
\emline{181.43}{10.89}{139}{196.43}{10.89}{140}
\emline{181.43}{0.89}{141}{196.43}{0.89}{142}
\put(186.43,5.89){\circle*{2.00}}
\put(191.43,5.89){\circle*{2.00}}
\put(29.43,22.89){\makebox(0,0)[cc]{$_6$}}
\put(34.43,22.89){\makebox(0,0)[cc]{$_4$}}
\put(39.43,22.89){\makebox(0,0)[cc]{$_2$}}
\put(29.43,17.89){\makebox(0,0)[cc]{$_7$}}
\put(34.43,17.89){\makebox(0,0)[cc]{$_5$}}
\put(39.43,17.89){\makebox(0,0)[cc]{$_3$}}
\put(44.43,22.89){\makebox(0,0)[cc]{$_0$}}
\put(44.43,17.89){\makebox(0,0)[cc]{$_1$}}
\put(49.43,22.89){\makebox(0,0)[cc]{$_6$}}
\put(54.43,22.89){\makebox(0,0)[cc]{$_4$}}
\put(59.43,22.89){\makebox(0,0)[cc]{$_2$}}
\put(49.43,17.89){\makebox(0,0)[cc]{$_7$}}
\put(54.43,17.89){\makebox(0,0)[cc]{$_5$}}
\put(59.43,17.89){\makebox(0,0)[cc]{$_3$}}
\put(64.43,22.89){\makebox(0,0)[cc]{$_0$}}
\put(64.43,17.89){\makebox(0,0)[cc]{$_1$}}
\put(69.43,22.89){\makebox(0,0)[cc]{$_6$}}
\put(74.43,22.89){\makebox(0,0)[cc]{$_4$}}
\put(79.43,22.89){\makebox(0,0)[cc]{$_2$}}
\put(69.43,17.89){\makebox(0,0)[cc]{$_7$}}
\put(74.43,17.89){\makebox(0,0)[cc]{$_5$}}
\put(79.43,17.89){\makebox(0,0)[cc]{$_3$}}
\put(84.43,22.89){\makebox(0,0)[cc]{$_0$}}
\put(84.43,17.89){\makebox(0,0)[cc]{$_1$}}
\put(89.43,22.89){\makebox(0,0)[cc]{$_6$}}
\put(94.43,22.89){\makebox(0,0)[cc]{$_4$}}
\put(99.43,22.89){\makebox(0,0)[cc]{$_2$}}
\put(89.43,17.89){\makebox(0,0)[cc]{$_7$}}
\put(94.43,17.89){\makebox(0,0)[cc]{$_5$}}
\put(99.43,17.89){\makebox(0,0)[cc]{$_3$}}
\put(104.43,17.89){\makebox(0,0)[cc]{$_1$}}
\put(109.43,17.89){\makebox(0,0)[cc]{$_0$}}
\put(114.43,17.89){\makebox(0,0)[cc]{$_6$}}
\put(119.43,17.89){\makebox(0,0)[cc]{$_4$}}
\put(124.43,17.89){\makebox(0,0)[cc]{$_2$}}
\put(114.43,12.89){\makebox(0,0)[cc]{$_7$}}
\put(119.43,12.89){\makebox(0,0)[cc]{$_5$}}
\put(124.43,12.89){\makebox(0,0)[cc]{$_3$}}
\put(129.43,12.89){\makebox(0,0)[cc]{$_1$}}
\put(134.43,12.89){\makebox(0,0)[cc]{$_0$}}
\put(139.43,12.89){\makebox(0,0)[cc]{$_6$}}
\put(144.43,12.89){\makebox(0,0)[cc]{$_4$}}
\put(149.43,12.89){\makebox(0,0)[cc]{$_2$}}
\put(139.43,7.89){\makebox(0,0)[cc]{$_7$}}
\put(144.43,7.89){\makebox(0,0)[cc]{$_5$}}
\put(149.43,7.89){\makebox(0,0)[cc]{$_3$}}
\put(154.43,7.89){\makebox(0,0)[cc]{$_1$}}
\put(159.43,7.89){\makebox(0,0)[cc]{$_0$}}
\put(164.43,7.89){\makebox(0,0)[cc]{$_6$}}
\put(169.43,7.89){\makebox(0,0)[cc]{$_4$}}
\put(174.43,7.89){\makebox(0,0)[cc]{$_2$}}
\put(164.43,2.89){\makebox(0,0)[cc]{$_7$}}
\put(169.43,2.89){\makebox(0,0)[cc]{$_5$}}
\put(174.43,2.89){\makebox(0,0)[cc]{$_3$}}
\put(179.43,7.89){\makebox(0,0)[cc]{$_0$}}
\put(179.43,2.89){\makebox(0,0)[cc]{$_1$}}
\put(184.43,7.89){\makebox(0,0)[cc]{$_6$}}
\put(189.43,7.89){\makebox(0,0)[cc]{$_4$}}
\put(194.43,7.89){\makebox(0,0)[cc]{$_2$}}
\put(184.43,2.89){\makebox(0,0)[cc]{$_7$}}
\put(189.43,2.89){\makebox(0,0)[cc]{$_5$}}
\put(194.43,2.89){\makebox(0,0)[cc]{$_3$}}
\put(199.43,2.89){\makebox(0,0)[cc]{$_1$}}
\put(204.43,2.89){\makebox(0,0)[cc]{$_0$}}
\emline{31.43}{25.89}{143}{31.43}{26.89}{144}
\emline{36.43}{25.89}{145}{36.43}{26.89}{146}
\emline{51.43}{25.89}{147}{51.43}{26.89}{148}
\emline{56.43}{25.89}{149}{56.43}{26.89}{150}
\emline{71.43}{25.89}{151}{71.43}{26.89}{152}
\emline{76.43}{25.89}{153}{76.43}{26.89}{154}
\emline{91.43}{25.89}{155}{91.43}{26.89}{156}
\emline{96.43}{25.89}{157}{96.43}{26.89}{158}
\emline{31.43}{14.89}{159}{31.43}{15.89}{160}
\emline{36.43}{14.89}{161}{36.43}{15.89}{162}
\emline{51.43}{14.89}{163}{51.43}{15.89}{164}
\emline{56.43}{14.89}{165}{56.43}{15.89}{166}
\emline{71.43}{14.89}{167}{71.43}{15.89}{168}
\emline{76.43}{14.89}{169}{76.43}{15.89}{170}
\emline{91.43}{14.89}{171}{91.43}{15.89}{172}
\emline{96.43}{14.89}{173}{96.43}{15.89}{174}
\emline{116.43}{20.89}{175}{116.43}{21.89}{176}
\emline{121.43}{20.89}{177}{121.43}{21.89}{178}
\emline{116.43}{9.89}{179}{116.43}{10.89}{180}
\emline{121.43}{9.89}{181}{121.43}{10.89}{182}
\emline{141.43}{15.89}{183}{141.43}{16.89}{184}
\emline{146.43}{15.89}{185}{146.43}{16.89}{186}
\emline{141.43}{4.89}{187}{141.43}{5.89}{188}
\emline{146.43}{4.89}{189}{146.43}{5.89}{190}
\emline{166.43}{10.89}{191}{166.43}{11.89}{192}
\emline{171.43}{10.89}{193}{171.43}{11.89}{194}
\emline{166.43}{-0.11}{195}{166.43}{0.89}{196}
\emline{171.43}{-0.11}{197}{171.43}{0.89}{198}
\emline{186.43}{10.89}{199}{186.43}{11.89}{200}
\emline{191.43}{10.89}{201}{191.43}{11.89}{202}
\emline{186.43}{-0.11}{203}{186.43}{0.89}{204}
\emline{191.43}{-0.11}{205}{191.43}{0.89}{206}
\emline{196.43}{10.89}{207}{196.43}{11.89}{208}
\emline{151.43}{15.89}{209}{151.43}{16.89}{210}
\emline{136.43}{4.89}{211}{136.43}{5.89}{212}
\emline{161.43}{-0.11}{213}{161.43}{0.89}{214}
\emline{206.43}{5.89}{215}{207.43}{5.89}{216}
\emline{206.43}{-0.11}{217}{206.43}{0.89}{218}
\emline{101.43}{25.89}{219}{101.43}{26.89}{220}
\emline{126.43}{20.89}{221}{126.43}{21.89}{222}
\emline{111.43}{9.89}{223}{111.43}{10.89}{224}
\emline{26.43}{25.89}{225}{25.43}{25.89}{226}
\emline{25.43}{20.89}{227}{26.43}{20.89}{228}
\emline{26.43}{14.89}{229}{26.43}{15.89}{230}
\end{picture}
\caption{The fundamental tile ${\mathcal T}(00011101)$}
\end{figure}
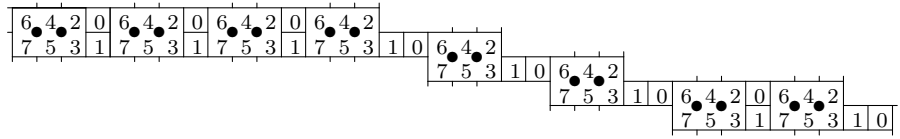

The unit-area squares
determined by the 4-cycles of ${\Lambda}$ in the Euclidean plane
are labeled for convenience as shown
in Figure 2, which establishes a fixed labeling of the realizations of the
pairs $(32,12)$ and $(32,21)$
in $\Psi(B)$. By extending the period $B$ to left and right and
applying $s$ and its powers, a larger picture of this labeling
representing $\Psi(B)$ can be seen
as a $32\times 32$ array $M=M(00011101)$ produced by
these labels of unit-area squares, whose first two rows
form the sub-array $M_1=M_1(00011101)=$
$$\begin{array}{c}
^{\,\,6\,\,4\,\,2\,\,0\,\,6\,\,4\,\,2\,\,0\,\,6\,\,4\,\,2\,\,0\,\,6\,\,4\,\,2\,\,7\,\,5\,\,3\,\,1\,\,0\,\,6\,\,4\,\,2\,\,7\,\,5\,\,3\,\,1\,\,7\,\,5\,\,3\,\,1\,\,0}
_{\,\,7\,\,5\,\,3\,\,1\,\,7\,\,5\,\,3\,\,1\,\,7\,\,5\,\,3\,\,1\,\,7\,\,5\,\,3\,\,1\,\,0\,\,6\,\,4\,\,2\,\,7\,\,5\,\,3\,\,1\,\,0\,\,6\,\,4\,\,2\,\,0\,\,6\,\,4\,\,2}
\end{array}$$
from which the whole
array can be obtained by means of the following induction step,
where $j=0,\ldots,15$:
Let $M_j$ be the sub-array of $M$ formed by the $(2j+1)$-th and $(2j+2)$-th rows of M.
If $M_j=(X,Y)$, where $Y$ is the rightmost $2\times 2$ sub-array of $M_j$, then
$M_{j+1}=(Y,X)$. This induction step is produced by the transformation $s$ and
what follows.

The union of the $32\times 32=1024$ 4-cycle-delimiting unit
squares whose labels are arranged precisely as in $M(B)$ form a
square $N(B)$ of side length 32 in the Euclidean plane. By
identifying the top and bottom sides of $N(B)$ as well as its left
and right sides, a (flat) torus ${\mathcal S}$ is obtained, (flat torus in the
sense of \cite{Costa}). In ${\mathcal S}$, the
projections of the original vertices and edges of ${\Lambda}$
form the Cartesian product $C_p\times C_p$ ($p=32$). In this $C_p\times C_p$,
the PTPC $S$ in ${\Lambda}$ determined by the period
$B=00011101$ projects onto a PTPC $\overline{S}$ in $C_p\times C_p$,
as required in the statement of the theorem. Also, $N(B)$ can be
furnished as a {\it cutout} of ${\mathcal S}$ in which $S$ or
$\overline{S}$ can be visualized, as is done down below in Figures
4 and 5 for the first two even- and odd-weight cases of periods
$B$ in their listings as indicated above, respectively.

Given a PDS $S$ in ${\Lambda}$, a labeling $f_S$ of ${\Lambda}$
with symbols in $\{0,1,2,3,4\}$ is defined as
follows: let 2 label each $v\in S$;
for each vertex $v$ of ${\Lambda}\setminus S$, let 0,1,3,4 label $v$ if
$v-(0,1),v-(1,0),v+(1,0),v+(0,1)$ belongs respectively to $S$.

The exemplified period $B=00011101$ contains all eight cases of
possible binary triples.
The vertex labels of $f_S$ associated to the room-ladder pairs in the image of the
middle entries of such triples through $\Psi$ are as shown in Figure 3,
(where also the unit-area squares determined by the 4-cycles of ${\Lambda}$
are labeled in their centers as indicated above).

\begin{figure}
\unitlength=0.70mm
\special{em:linewidth 0.4pt}
\linethickness{0.4pt}
\begin{picture}(171.04,63.82)
\put(48.04,58.82){\circle{2.00}}
\put(58.04,58.82){\circle{2.00}}
\put(68.04,58.82){\circle{2.00}}
\put(78.04,58.82){\circle{2.00}}
\put(78.04,48.82){\circle{2.00}}
\put(78.04,38.82){\circle{2.00}}
\put(68.04,38.82){\circle{2.00}}
\put(58.04,38.82){\circle{2.00}}
\put(48.04,38.82){\circle{2.00}}
\put(48.04,48.82){\circle{2.00}}
\put(58.04,48.82){\circle*{2.00}}
\put(68.04,48.82){\circle*{2.00}}
\put(88.04,58.82){\circle{2.00}}
\put(88.04,48.82){\circle{2.00}}
\put(88.04,38.82){\circle{2.00}}
\emline{49.04}{58.82}{1}{57.04}{58.82}{2}
\emline{59.04}{58.82}{3}{67.04}{58.82}{4}
\emline{69.04}{58.82}{5}{77.04}{58.82}{6}
\emline{79.04}{58.82}{7}{87.04}{58.82}{8}
\emline{88.04}{57.82}{9}{88.04}{49.82}{10}
\emline{88.04}{47.82}{11}{88.04}{39.82}{12}
\emline{87.04}{38.82}{13}{79.04}{38.82}{14}
\emline{77.04}{38.82}{15}{69.04}{38.82}{16}
\emline{67.04}{38.82}{17}{59.04}{38.82}{18}
\emline{57.04}{38.82}{19}{49.04}{38.82}{20}
\emline{48.04}{39.82}{21}{48.04}{47.82}{22}
\emline{48.04}{49.82}{23}{48.04}{57.82}{24}
\emline{78.04}{57.82}{25}{78.04}{49.82}{26}
\emline{78.04}{47.82}{27}{78.04}{39.82}{28}
\put(51.04,55.82){\makebox(0,0)[cc]{$_4$}}
\put(61.04,55.82){\makebox(0,0)[cc]{$_0$}}
\put(71.04,55.82){\makebox(0,0)[cc]{$_0$}}
\put(81.04,55.82){\makebox(0,0)[cc]{$_4$}}
\put(91.04,55.82){\makebox(0,0)[cc]{$_4$}}
\put(51.04,45.82){\makebox(0,0)[cc]{$_1$}}
\put(61.04,45.82){\makebox(0,0)[cc]{$_2$}}
\put(71.04,45.82){\makebox(0,0)[cc]{$_2$}}
\put(81.04,45.82){\makebox(0,0)[cc]{$_3$}}
\put(91.04,45.82){\makebox(0,0)[cc]{$_1$}}
\put(51.04,35.82){\makebox(0,0)[cc]{$_0$}}
\put(61.04,35.82){\makebox(0,0)[cc]{$_4$}}
\put(71.04,35.82){\makebox(0,0)[cc]{$_4$}}
\put(81.04,35.82){\makebox(0,0)[cc]{$_0$}}
\put(91.04,35.82){\makebox(0,0)[cc]{$_0$}}
\put(68.04,63.82){\makebox(0,0)[cc]{$_{000,\,\,001:}$}}
\put(54.04,52.82){\makebox(0,0)[cc]{$_6$}}
\put(64.04,52.82){\makebox(0,0)[cc]{$_4$}}
\put(74.04,52.82){\makebox(0,0)[cc]{$_2$}}
\put(84.04,52.82){\makebox(0,0)[cc]{$_0$}}
\put(54.04,42.82){\makebox(0,0)[cc]{$_7$}}
\put(64.04,42.82){\makebox(0,0)[cc]{$_5$}}
\put(74.04,42.82){\makebox(0,0)[cc]{$_3$}}
\put(84.04,42.82){\makebox(0,0)[cc]{$_1$}}
\put(118.04,58.82){\circle{2.00}}
\put(128.04,58.82){\circle{2.00}}
\put(138.04,58.82){\circle{2.00}}
\put(148.04,58.82){\circle{2.00}}
\put(148.04,48.82){\circle{2.00}}
\put(148.04,38.82){\circle{2.00}}
\put(138.04,38.82){\circle{2.00}}
\put(128.04,38.82){\circle{2.00}}
\put(118.04,38.82){\circle{2.00}}
\put(118.04,48.82){\circle{2.00}}
\put(128.04,48.82){\circle*{2.00}}
\put(138.04,48.82){\circle*{2.00}}
\put(158.04,58.82){\circle{2.00}}
\put(158.04,48.82){\circle{2.00}}
\put(158.04,38.82){\circle{2.00}}
\emline{119.04}{58.82}{29}{127.04}{58.82}{30}
\emline{129.04}{58.82}{31}{137.04}{58.82}{32}
\emline{139.04}{58.82}{33}{147.04}{58.82}{34}
\emline{149.04}{58.82}{35}{157.04}{58.82}{36}
\emline{158.04}{57.82}{37}{158.04}{49.82}{38}
\emline{158.04}{47.82}{39}{158.04}{39.82}{40}
\emline{157.04}{38.82}{41}{149.04}{38.82}{42}
\emline{147.04}{38.82}{43}{139.04}{38.82}{44}
\emline{137.04}{38.82}{45}{129.04}{38.82}{46}
\emline{127.04}{38.82}{47}{119.04}{38.82}{48}
\emline{118.04}{39.82}{49}{118.04}{47.82}{50}
\emline{118.04}{49.82}{51}{118.04}{57.82}{52}
\emline{148.04}{57.82}{53}{148.04}{49.82}{54}
\emline{148.04}{47.82}{55}{148.04}{39.82}{56}
\put(121.04,55.82){\makebox(0,0)[cc]{$_4$}}
\put(131.04,55.82){\makebox(0,0)[cc]{$_0$}}
\put(141.04,55.82){\makebox(0,0)[cc]{$_0$}}
\put(151.04,55.82){\makebox(0,0)[cc]{$_4$}}
\put(161.04,55.82){\makebox(0,0)[cc]{$_4$}}
\put(121.04,45.82){\makebox(0,0)[cc]{$_1$}}
\put(131.04,45.82){\makebox(0,0)[cc]{$_2$}}
\put(141.04,45.82){\makebox(0,0)[cc]{$_2$}}
\put(151.04,45.82){\makebox(0,0)[cc]{$_3$}}
\put(161.04,45.82){\makebox(0,0)[cc]{$_1$}}
\put(121.04,35.82){\makebox(0,0)[cc]{$_3$}}
\put(131.04,35.82){\makebox(0,0)[cc]{$_4$}}
\put(141.04,35.82){\makebox(0,0)[cc]{$_4$}}
\put(151.04,35.82){\makebox(0,0)[cc]{$_0$}}
\put(161.04,35.82){\makebox(0,0)[cc]{$_0$}}
\put(138.04,63.82){\makebox(0,0)[cc]{$_{100,\,\,101:}$}}
\put(124.04,52.82){\makebox(0,0)[cc]{$_6$}}
\put(134.04,52.82){\makebox(0,0)[cc]{$_4$}}
\put(144.04,52.82){\makebox(0,0)[cc]{$_2$}}
\put(154.04,52.82){\makebox(0,0)[cc]{$_0$}}
\put(124.04,42.82){\makebox(0,0)[cc]{$_7$}}
\put(134.04,42.82){\makebox(0,0)[cc]{$_5$}}
\put(144.04,42.82){\makebox(0,0)[cc]{$_3$}}
\put(154.04,42.82){\makebox(0,0)[cc]{$_1$}}
\put(48.04,23.82){\circle{2.00}}
\put(58.04,23.82){\circle{2.00}}
\put(68.04,23.82){\circle{2.00}}
\put(78.04,23.82){\circle{2.00}}
\put(78.04,13.82){\circle{2.00}}
\put(78.04,3.82){\circle{2.00}}
\put(68.04,3.82){\circle{2.00}}
\put(58.04,3.82){\circle{2.00}}
\put(48.04,3.82){\circle{2.00}}
\put(48.04,13.82){\circle{2.00}}
\put(58.04,13.82){\circle*{2.00}}
\put(68.04,13.82){\circle*{2.00}}
\put(88.04,13.82){\circle{2.00}}
\put(88.04,3.82){\circle{2.00}}
\emline{49.04}{23.82}{57}{57.04}{23.82}{58}
\emline{59.04}{23.82}{59}{67.04}{23.82}{60}
\emline{69.04}{23.82}{61}{77.04}{23.82}{62}
\emline{88.04}{12.82}{63}{88.04}{4.82}{64}
\emline{87.04}{3.82}{65}{79.04}{3.82}{66}
\emline{77.04}{3.82}{67}{69.04}{3.82}{68}
\emline{67.04}{3.82}{69}{59.04}{3.82}{70}
\emline{57.04}{3.82}{71}{49.04}{3.82}{72}
\emline{48.04}{4.82}{73}{48.04}{12.82}{74}
\emline{48.04}{14.82}{75}{48.04}{22.82}{76}
\emline{78.04}{22.82}{77}{78.04}{14.82}{78}
\emline{78.04}{12.82}{79}{78.04}{4.82}{80}
\put(51.04,20.82){\makebox(0,0)[cc]{$_4$}}
\put(61.04,20.82){\makebox(0,0)[cc]{$_0$}}
\put(71.04,20.82){\makebox(0,0)[cc]{$_0$}}
\put(81.04,20.82){\makebox(0,0)[cc]{$_1$}}
\put(51.04,10.82){\makebox(0,0)[cc]{$_1$}}
\put(61.04,10.82){\makebox(0,0)[cc]{$_2$}}
\put(71.04,10.82){\makebox(0,0)[cc]{$_2$}}
\put(81.04,10.82){\makebox(0,0)[cc]{$_3$}}
\put(91.04,10.82){\makebox(0,0)[cc]{$_4$}}
\put(51.04,0.82){\makebox(0,0)[cc]{$_0$}}
\put(61.04,0.82){\makebox(0,0)[cc]{$_4$}}
\put(71.04,0.82){\makebox(0,0)[cc]{$_4$}}
\put(81.04,0.82){\makebox(0,0)[cc]{$_0$}}
\put(91.04,0.82){\makebox(0,0)[cc]{$_0$}}
\put(68.04,28.82){\makebox(0,0)[cc]{$_{010,\,\,011:}$}}
\put(54.04,17.82){\makebox(0,0)[cc]{$_6$}}
\put(64.04,17.82){\makebox(0,0)[cc]{$_4$}}
\put(74.04,17.82){\makebox(0,0)[cc]{$_2$}}
\put(54.04,7.82){\makebox(0,0)[cc]{$_7$}}
\put(64.04,7.82){\makebox(0,0)[cc]{$_5$}}
\put(74.04,7.82){\makebox(0,0)[cc]{$_3$}}
\put(84.04,7.82){\makebox(0,0)[cc]{$_1$}}
\put(118.04,23.82){\circle{2.00}}
\put(128.04,23.82){\circle{2.00}}
\put(138.04,23.82){\circle{2.00}}
\put(148.04,23.82){\circle{2.00}}
\put(148.04,13.82){\circle{2.00}}
\put(148.04,3.82){\circle{2.00}}
\put(138.04,3.82){\circle{2.00}}
\put(128.04,3.82){\circle{2.00}}
\put(118.04,3.82){\circle{2.00}}
\put(118.04,13.82){\circle{2.00}}
\put(128.04,13.82){\circle*{2.00}}
\put(138.04,13.82){\circle*{2.00}}
\put(158.04,13.82){\circle{2.00}}
\put(158.04,3.82){\circle{2.00}}
\emline{119.04}{23.82}{81}{127.04}{23.82}{82}
\emline{129.04}{23.82}{83}{137.04}{23.82}{84}
\emline{139.04}{23.82}{85}{147.04}{23.82}{86}
\emline{158.04}{12.82}{87}{158.04}{4.82}{88}
\emline{157.04}{3.82}{89}{149.04}{3.82}{90}
\emline{147.04}{3.82}{91}{139.04}{3.82}{92}
\emline{137.04}{3.82}{93}{129.04}{3.82}{94}
\emline{127.04}{3.82}{95}{119.04}{3.82}{96}
\emline{118.04}{4.82}{97}{118.04}{12.82}{98}
\emline{118.04}{14.82}{99}{118.04}{22.82}{100}
\emline{148.04}{22.82}{101}{148.04}{14.82}{102}
\emline{148.04}{12.82}{103}{148.04}{4.82}{104}
\put(121.04,20.82){\makebox(0,0)[cc]{$_4$}}
\put(131.04,20.82){\makebox(0,0)[cc]{$_0$}}
\put(141.04,20.82){\makebox(0,0)[cc]{$_0$}}
\put(151.04,20.82){\makebox(0,0)[cc]{$_1$}}
\put(121.04,10.82){\makebox(0,0)[cc]{$_1$}}
\put(131.04,10.82){\makebox(0,0)[cc]{$_2$}}
\put(141.04,10.82){\makebox(0,0)[cc]{$_2$}}
\put(151.04,10.82){\makebox(0,0)[cc]{$_3$}}
\put(161.04,10.82){\makebox(0,0)[cc]{$_4$}}
\put(121.04,0.82){\makebox(0,0)[cc]{$_3$}}
\put(131.04,0.82){\makebox(0,0)[cc]{$_4$}}
\put(141.04,0.82){\makebox(0,0)[cc]{$_4$}}
\put(151.04,0.82){\makebox(0,0)[cc]{$_0$}}
\put(161.04,0.82){\makebox(0,0)[cc]{$_0$}}
\put(138.04,28.82){\makebox(0,0)[cc]{$_{110,\,\,111:}$}}
\put(124.04,17.82){\makebox(0,0)[cc]{$_6$}}
\put(134.04,17.82){\makebox(0,0)[cc]{$_4$}}
\put(144.04,17.82){\makebox(0,0)[cc]{$_2$}}
\put(124.04,7.82){\makebox(0,0)[cc]{$_7$}}
\put(134.04,7.82){\makebox(0,0)[cc]{$_5$}}
\put(144.04,7.82){\makebox(0,0)[cc]{$_3$}}
\put(154.04,7.82){\makebox(0,0)[cc]{$_1$}}
\emline{79.04}{48.82}{105}{87.04}{48.82}{106}
\emline{149.04}{48.82}{107}{157.04}{48.82}{108}
\emline{157.04}{13.82}{109}{149.04}{13.82}{110}
\emline{87.04}{13.82}{111}{79.04}{13.82}{112}
\put(98.04,13.82){\circle{2.00}}
\put(98.04,3.82){\circle{2.00}}
\emline{89.04}{13.82}{113}{97.04}{13.82}{114}
\emline{98.04}{12.82}{115}{98.04}{4.82}{116}
\put(101.04,10.82){\makebox(0,0)[cc]{$_4$}}
\put(94.04,7.82){\makebox(0,0)[cc]{$_0$}}
\emline{97.04}{3.82}{117}{89.04}{3.82}{118}
\put(168.04,13.82){\circle{2.00}}
\put(168.04,3.82){\circle{2.00}}
\emline{159.04}{13.82}{119}{167.04}{13.82}{120}
\emline{168.04}{12.82}{121}{168.04}{4.82}{122}
\put(171.04,10.82){\makebox(0,0)[cc]{$_4$}}
\put(171.04,0.82){\makebox(0,0)[cc]{$_1$}}
\put(164.04,7.82){\makebox(0,0)[cc]{$_0$}}
\emline{167.04}{3.82}{123}{159.04}{3.82}{124}
\put(101.04,0.82){\makebox(0,0)[cc]{$_1$}}
\end{picture}
\caption{Labelled room-label pairs for middle entries of binary triples}
\end{figure}
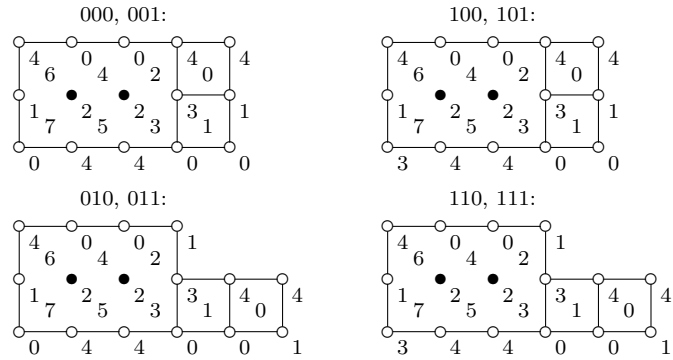

In general, a fundamental tile ${\mathcal T}(B)$
is determined, for any period $B$ of a periodic
doubly infinite binary sequence by replacing each entry $b_i=0$ of $B$ by a sub-tile
$\tau_0$ as in the two top cases of Figure 3, formed by a
room-ladder pair with PDS-sub-array $(32,12)$, and by replacing
each entry $b_i=1$ of $B$
by a sub-tile $\tau_1$ as in the two bottom cases of Figure 3, formed by a
room-ladder pair with PDS-sub-array $(32,21)$. If $b_i=0$, then
the unit-area squares with labels 0 and 1 in $\tau_0$ are adjacent to their right to
unit-area squares with respective labels 7 and 6, belonging to a sub-tile of either
type. if $b_i=1$, then the unit-area square with label 0 is adjacent to its right to a
unit-area square with label 7, that is
at the upper-left corner of another sub-tile. This composed periodically (for a fixed
period B) and translated by the powers of $s$ yields $\Psi(B)$.

Given a period $B$ of a periodic doubly infinite binary sequence,
let $B_1$ be the subsequence of $B$ whose last entry is the first, or leftmost,
unit entry of $B$. Inductively for $1<i\leq k$,
let $B_i$ be the subsequence of
$B\setminus(\cup_{j=1}^{i-1} B_j)$ whose last entry is the $i$-th unit entry of $B$,
where $B\setminus(\cup_{j=1}^{i-1} B_j)$ is the subsequence of $B$ obtained by deleting
the successive concatenation of $B_1,B_2,\ldots,B_{i-1}$.
The sequence $\xi(B_1)$
of labels of unit-area squares in the first, or top, row of ${\mathcal T}$
is obtained by successively replacing each null entry of $B_1$ by a subsequence 6420,
and the final unit entry of $B_1$
by a subsequence 642. There is an accompaniment of this $\xi(B_1)$ by a sequence $\eta(B_1)$
in the second row of such labels in ${\mathcal T}(B)$: below each subsequence
6420 in $\xi(B_1)$,
corresponds a subsequence 7531; below the final subsequence 642 of $\xi(B_1)$
corresponds a subsequence 75310, extended two
extra positions to the right of the final entry of $\xi(B_1)$.
From the subsequent position to the right of this final entry of $\xi(B_1)$,
${\mathcal T}$ contains $\xi(B_2)$, and below it, $\eta(B_2)$, and so on, down to
$\xi(B_k)$ and $\eta(B_k)$.

For $B=00011101$,
the first row of $M_1(B)$ is formed by the concatenation of the first, third and fifth
rows of ${\mathcal B}$, namely $\xi(B_1)=\xi(0001)=(6420)^3642$, 
$\eta(B_2)\xi(B_3)=\eta(1)\xi(1)=75310642$ and $\eta(B_4)=\eta(01)=(7531)^20$, 
respectively;
the second row of $M_1(B)$ is formed by the concatenation of the second and the fourth
rows of ${\mathcal B}$, namely
$\eta(B_1)\xi(B_2)=\eta(0001)\xi(1)=(7531)^4 0642$ and
$\eta(B_3)\xi(B_4)=\eta(1)\xi(01)=7531(0642)^2$, 
respectively.

Similarly for any other period $B$ of even weight $k$: the first
(second) row of $M_1(B)$ is the concatenation of the odd (even) rows
of ${\mathcal T}(B)$ in their increasing order, which is then
repeated periodically in ${\Lambda}$ by way of horizontal
concatenations. Observe that the rows of ${\mathcal T}(B)$ have
successively the following label sequences, or concatenations of label sequences:
$$\xi(B_1),\eta(B_1)\xi(B_2),\ldots,\eta(B_i)\xi(B_{i+1}),\ldots,
\eta(B_{k-1})\xi(B_k),\eta(B_k),$$
from which we must select the concatenation of the odd- (even-)
positioned member label subsequences as the first (second) row
of $M_1(B)$, whose total length is then seen to be $4k$, ($4k$).

For each period $B$ of odd weight $k$,
the first (second) row of $M_1(B)$ is the concatenation of the odd (even) rows of
${\mathcal T}(B)$ followed by the concatenation of the even (odd) rows of
${\mathcal T}(B)$, in their increasing order in both cases.

The total length of the concatenations in either the first of the second
row is $p=4n$ if the weight $k$ of $B$ is even and is $p=8n$ if $k$ is odd.
This implies the statement of the theorem.
\qfd

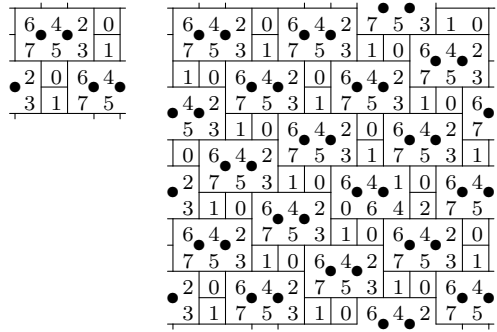
\begin{figure}
\unitlength=0.70mm
\special{em:linewidth 0.4pt}
\linethickness{0.4pt}
\begin{picture}(153.42,62.00)
\put(65.42,58.00){\makebox(0,0)[cc]{$_6$}}
\put(70.42,58.00){\makebox(0,0)[cc]{$_4$}}
\put(75.42,58.00){\makebox(0,0)[cc]{$_2$}}
\put(80.42,58.00){\makebox(0,0)[cc]{$_0$}}
\put(65.42,53.00){\makebox(0,0)[cc]{$_7$}}
\put(70.42,53.00){\makebox(0,0)[cc]{$_5$}}
\put(75.42,53.00){\makebox(0,0)[cc]{$_3$}}
\put(80.42,53.00){\makebox(0,0)[cc]{$_1$}}
\put(65.42,48.00){\makebox(0,0)[cc]{$_2$}}
\put(70.42,48.00){\makebox(0,0)[cc]{$_0$}}
\put(75.42,48.00){\makebox(0,0)[cc]{$_6$}}
\put(80.42,48.00){\makebox(0,0)[cc]{$_4$}}
\put(65.42,43.00){\makebox(0,0)[cc]{$_3$}}
\put(70.42,43.00){\makebox(0,0)[cc]{$_1$}}
\put(75.42,43.00){\makebox(0,0)[cc]{$_7$}}
\put(80.42,43.00){\makebox(0,0)[cc]{$_5$}}
\emline{77.42}{61.00}{1}{77.42}{51.00}{2}
\emline{82.42}{51.00}{3}{82.42}{61.00}{4}
\emline{82.42}{61.00}{5}{62.42}{61.00}{6}
\emline{62.42}{61.00}{7}{62.42}{51.00}{8}
\emline{62.42}{51.00}{9}{82.42}{51.00}{10}
\emline{67.42}{51.00}{11}{67.42}{41.00}{12}
\emline{62.42}{41.00}{13}{82.42}{41.00}{14}
\emline{72.42}{41.00}{15}{72.42}{51.00}{16}
\emline{72.42}{46.00}{17}{67.42}{46.00}{18}
\emline{77.42}{56.00}{19}{82.42}{56.00}{20}
\put(67.42,56.00){\circle*{2.00}}
\put(72.42,56.00){\circle*{2.00}}
\put(77.42,46.00){\circle*{2.00}}
\put(82.42,46.00){\circle*{2.00}}
\put(62.42,46.00){\circle*{2.00}}
\put(95.42,58.00){\makebox(0,0)[cc]{$_6$}}
\put(100.42,58.00){\makebox(0,0)[cc]{$_4$}}
\put(105.42,58.00){\makebox(0,0)[cc]{$_2$}}
\put(110.42,58.00){\makebox(0,0)[cc]{$_0$}}
\put(95.42,53.00){\makebox(0,0)[cc]{$_7$}}
\put(100.42,53.00){\makebox(0,0)[cc]{$_5$}}
\put(105.42,53.00){\makebox(0,0)[cc]{$_3$}}
\put(110.42,53.00){\makebox(0,0)[cc]{$_1$}}
\put(95.42,48.00){\makebox(0,0)[cc]{$_1$}}
\put(100.42,48.00){\makebox(0,0)[cc]{$_0$}}
\put(105.42,48.00){\makebox(0,0)[cc]{$_6$}}
\put(110.42,48.00){\makebox(0,0)[cc]{$_4$}}
\put(95.42,43.00){\makebox(0,0)[cc]{$_4$}}
\put(100.42,43.00){\makebox(0,0)[cc]{$_2$}}
\put(105.42,43.00){\makebox(0,0)[cc]{$_7$}}
\put(110.42,43.00){\makebox(0,0)[cc]{$_5$}}
\put(115.42,58.00){\makebox(0,0)[cc]{$_6$}}
\put(120.42,58.00){\makebox(0,0)[cc]{$_4$}}
\put(125.42,58.00){\makebox(0,0)[cc]{$_2$}}
\put(130.42,58.00){\makebox(0,0)[cc]{$_7$}}
\put(115.42,53.00){\makebox(0,0)[cc]{$_7$}}
\put(120.42,53.00){\makebox(0,0)[cc]{$_5$}}
\put(125.42,53.00){\makebox(0,0)[cc]{$_3$}}
\put(130.42,53.00){\makebox(0,0)[cc]{$_1$}}
\put(115.42,48.00){\makebox(0,0)[cc]{$_2$}}
\put(120.42,48.00){\makebox(0,0)[cc]{$_0$}}
\put(125.42,48.00){\makebox(0,0)[cc]{$_6$}}
\put(130.42,48.00){\makebox(0,0)[cc]{$_4$}}
\put(115.42,43.00){\makebox(0,0)[cc]{$_3$}}
\put(120.42,43.00){\makebox(0,0)[cc]{$_1$}}
\put(125.42,43.00){\makebox(0,0)[cc]{$_7$}}
\put(130.42,43.00){\makebox(0,0)[cc]{$_5$}}
\put(95.42,38.00){\makebox(0,0)[cc]{$_5$}}
\put(100.42,38.00){\makebox(0,0)[cc]{$_3$}}
\put(105.42,38.00){\makebox(0,0)[cc]{$_1$}}
\put(110.42,38.00){\makebox(0,0)[cc]{$_0$}}
\put(95.42,33.00){\makebox(0,0)[cc]{$_0$}}
\put(100.42,33.00){\makebox(0,0)[cc]{$_6$}}
\put(105.42,33.00){\makebox(0,0)[cc]{$_4$}}
\put(110.42,33.00){\makebox(0,0)[cc]{$_2$}}
\put(95.42,28.00){\makebox(0,0)[cc]{$_2$}}
\put(100.42,28.00){\makebox(0,0)[cc]{$_7$}}
\put(105.42,28.00){\makebox(0,0)[cc]{$_5$}}
\put(110.42,28.00){\makebox(0,0)[cc]{$_3$}}
\put(95.42,23.00){\makebox(0,0)[cc]{$_3$}}
\put(100.42,23.00){\makebox(0,0)[cc]{$_1$}}
\put(105.42,23.00){\makebox(0,0)[cc]{$_0$}}
\put(110.42,23.00){\makebox(0,0)[cc]{$_6$}}
\put(115.42,38.00){\makebox(0,0)[cc]{$_6$}}
\put(120.42,38.00){\makebox(0,0)[cc]{$_4$}}
\put(125.42,38.00){\makebox(0,0)[cc]{$_2$}}
\put(130.42,38.00){\makebox(0,0)[cc]{$_0$}}
\put(115.42,33.00){\makebox(0,0)[cc]{$_7$}}
\put(120.42,33.00){\makebox(0,0)[cc]{$_5$}}
\put(125.42,33.00){\makebox(0,0)[cc]{$_3$}}
\put(130.42,33.00){\makebox(0,0)[cc]{$_1$}}
\put(115.42,28.00){\makebox(0,0)[cc]{$_1$}}
\put(120.42,28.00){\makebox(0,0)[cc]{$_0$}}
\put(125.42,28.00){\makebox(0,0)[cc]{$_6$}}
\put(130.42,28.00){\makebox(0,0)[cc]{$_4$}}
\put(115.42,23.00){\makebox(0,0)[cc]{$_4$}}
\put(120.42,23.00){\makebox(0,0)[cc]{$_2$}}
\put(125.42,23.00){\makebox(0,0)[cc]{$_0$}}
\put(130.42,23.00){\makebox(0,0)[cc]{$_6$}}
\put(135.42,58.00){\makebox(0,0)[cc]{$_5$}}
\put(140.42,58.00){\makebox(0,0)[cc]{$_3$}}
\put(145.42,58.00){\makebox(0,0)[cc]{$_1$}}
\put(150.42,58.00){\makebox(0,0)[cc]{$_0$}}
\put(135.42,53.00){\makebox(0,0)[cc]{$_0$}}
\put(140.42,53.00){\makebox(0,0)[cc]{$_6$}}
\put(145.42,53.00){\makebox(0,0)[cc]{$_4$}}
\put(150.42,53.00){\makebox(0,0)[cc]{$_2$}}
\put(135.42,48.00){\makebox(0,0)[cc]{$_2$}}
\put(140.42,48.00){\makebox(0,0)[cc]{$_7$}}
\put(145.42,48.00){\makebox(0,0)[cc]{$_5$}}
\put(150.42,48.00){\makebox(0,0)[cc]{$_3$}}
\put(135.42,43.00){\makebox(0,0)[cc]{$_3$}}
\put(140.42,43.00){\makebox(0,0)[cc]{$_1$}}
\put(145.42,43.00){\makebox(0,0)[cc]{$_0$}}
\put(150.42,43.00){\makebox(0,0)[cc]{$_6$}}
\put(135.42,38.00){\makebox(0,0)[cc]{$_6$}}
\put(140.42,38.00){\makebox(0,0)[cc]{$_4$}}
\put(145.42,38.00){\makebox(0,0)[cc]{$_2$}}
\put(150.42,38.00){\makebox(0,0)[cc]{$_7$}}
\put(135.42,33.00){\makebox(0,0)[cc]{$_7$}}
\put(140.42,33.00){\makebox(0,0)[cc]{$_5$}}
\put(145.42,33.00){\makebox(0,0)[cc]{$_3$}}
\put(150.42,33.00){\makebox(0,0)[cc]{$_1$}}
\put(135.42,28.00){\makebox(0,0)[cc]{$_1$}}
\put(140.42,28.00){\makebox(0,0)[cc]{$_0$}}
\put(145.42,28.00){\makebox(0,0)[cc]{$_6$}}
\put(150.42,28.00){\makebox(0,0)[cc]{$_4$}}
\put(135.42,23.00){\makebox(0,0)[cc]{$_4$}}
\put(140.42,23.00){\makebox(0,0)[cc]{$_2$}}
\put(145.42,23.00){\makebox(0,0)[cc]{$_7$}}
\put(150.42,23.00){\makebox(0,0)[cc]{$_5$}}
\put(95.42,18.00){\makebox(0,0)[cc]{$_6$}}
\put(100.42,18.00){\makebox(0,0)[cc]{$_4$}}
\put(105.42,18.00){\makebox(0,0)[cc]{$_2$}}
\put(110.42,18.00){\makebox(0,0)[cc]{$_7$}}
\put(95.42,13.00){\makebox(0,0)[cc]{$_7$}}
\put(100.42,13.00){\makebox(0,0)[cc]{$_5$}}
\put(105.42,13.00){\makebox(0,0)[cc]{$_3$}}
\put(110.42,13.00){\makebox(0,0)[cc]{$_1$}}
\put(95.42,8.00){\makebox(0,0)[cc]{$_2$}}
\put(100.42,8.00){\makebox(0,0)[cc]{$_0$}}
\put(105.42,8.00){\makebox(0,0)[cc]{$_6$}}
\put(110.42,8.00){\makebox(0,0)[cc]{$_4$}}
\put(95.42,3.00){\makebox(0,0)[cc]{$_3$}}
\put(100.42,3.00){\makebox(0,0)[cc]{$_1$}}
\put(105.42,3.00){\makebox(0,0)[cc]{$_7$}}
\put(110.42,3.00){\makebox(0,0)[cc]{$_5$}}
\put(115.42,18.00){\makebox(0,0)[cc]{$_5$}}
\put(120.42,18.00){\makebox(0,0)[cc]{$_3$}}
\put(125.42,18.00){\makebox(0,0)[cc]{$_1$}}
\put(130.42,18.00){\makebox(0,0)[cc]{$_0$}}
\put(115.42,13.00){\makebox(0,0)[cc]{$_0$}}
\put(120.42,13.00){\makebox(0,0)[cc]{$_6$}}
\put(125.42,13.00){\makebox(0,0)[cc]{$_4$}}
\put(130.42,13.00){\makebox(0,0)[cc]{$_2$}}
\put(115.42,8.00){\makebox(0,0)[cc]{$_2$}}
\put(120.42,8.00){\makebox(0,0)[cc]{$_7$}}
\put(125.42,8.00){\makebox(0,0)[cc]{$_5$}}
\put(130.42,8.00){\makebox(0,0)[cc]{$_3$}}
\put(115.42,3.00){\makebox(0,0)[cc]{$_3$}}
\put(120.42,3.00){\makebox(0,0)[cc]{$_1$}}
\put(125.42,3.00){\makebox(0,0)[cc]{$_0$}}
\put(130.42,3.00){\makebox(0,0)[cc]{$_6$}}
\put(135.42,18.00){\makebox(0,0)[cc]{$_6$}}
\put(140.42,18.00){\makebox(0,0)[cc]{$_4$}}
\put(145.42,18.00){\makebox(0,0)[cc]{$_2$}}
\put(150.42,18.00){\makebox(0,0)[cc]{$_0$}}
\put(135.42,13.00){\makebox(0,0)[cc]{$_7$}}
\put(140.42,13.00){\makebox(0,0)[cc]{$_5$}}
\put(145.42,13.00){\makebox(0,0)[cc]{$_3$}}
\put(150.42,13.00){\makebox(0,0)[cc]{$_1$}}
\put(135.42,8.00){\makebox(0,0)[cc]{$_1$}}
\put(140.42,8.00){\makebox(0,0)[cc]{$_0$}}
\put(145.42,8.00){\makebox(0,0)[cc]{$_6$}}
\put(150.42,8.00){\makebox(0,0)[cc]{$_4$}}
\put(135.42,3.00){\makebox(0,0)[cc]{$_4$}}
\put(140.42,3.00){\makebox(0,0)[cc]{$_2$}}
\put(145.42,3.00){\makebox(0,0)[cc]{$_7$}}
\put(150.42,3.00){\makebox(0,0)[cc]{$_5$}}
\emline{92.42}{61.00}{21}{127.42}{61.00}{22}
\emline{127.42}{61.00}{23}{127.42}{51.00}{24}
\emline{127.42}{51.00}{25}{92.42}{51.00}{26}
\emline{92.42}{51.00}{27}{92.42}{61.00}{28}
\emline{107.42}{61.00}{29}{107.42}{51.00}{30}
\emline{107.42}{56.00}{31}{112.42}{56.00}{32}
\emline{112.42}{61.00}{33}{112.42}{51.00}{34}
\emline{127.42}{51.00}{35}{137.42}{51.00}{36}
\emline{137.42}{51.00}{37}{137.42}{56.00}{38}
\emline{137.42}{56.00}{39}{127.42}{56.00}{40}
\emline{132.42}{56.00}{41}{132.42}{51.00}{42}
\emline{137.42}{51.00}{43}{137.42}{46.00}{44}
\emline{137.42}{46.00}{45}{152.42}{46.00}{46}
\emline{152.42}{46.00}{47}{152.42}{56.00}{48}
\emline{152.42}{56.00}{49}{137.42}{56.00}{50}
\emline{142.42}{56.00}{51}{142.42}{61.00}{52}
\emline{142.42}{61.00}{53}{152.42}{61.00}{54}
\emline{152.42}{61.00}{55}{152.42}{56.00}{56}
\put(147.42,56.00){\rule{0.00\unitlength}{5.00\unitlength}}
\emline{92.42}{51.00}{57}{92.42}{46.00}{58}
\emline{92.42}{46.00}{59}{102.42}{46.00}{60}
\emline{102.42}{46.00}{61}{102.42}{51.00}{62}
\emline{97.42}{51.00}{63}{97.42}{46.00}{64}
\emline{102.42}{46.00}{65}{102.42}{41.00}{66}
\emline{102.42}{41.00}{67}{117.42}{41.00}{68}
\emline{117.42}{41.00}{69}{117.42}{51.00}{70}
\emline{117.42}{41.00}{71}{122.42}{41.00}{72}
\emline{122.42}{41.00}{73}{122.42}{51.00}{74}
\emline{117.42}{46.00}{75}{122.42}{46.00}{76}
\emline{122.42}{41.00}{77}{147.42}{41.00}{78}
\emline{147.42}{41.00}{79}{147.42}{46.00}{80}
\emline{142.42}{46.00}{81}{142.42}{41.00}{82}
\emline{137.42}{41.00}{83}{137.42}{46.00}{84}
\emline{147.42}{41.00}{85}{147.42}{36.00}{86}
\emline{147.42}{36.00}{87}{152.42}{36.00}{88}
\emline{92.42}{36.00}{89}{102.42}{36.00}{90}
\emline{102.42}{36.00}{91}{102.42}{41.00}{92}
\emline{102.42}{36.00}{93}{112.42}{36.00}{94}
\emline{112.42}{36.00}{95}{112.42}{41.00}{96}
\emline{107.42}{41.00}{97}{107.42}{36.00}{98}
\emline{112.42}{36.00}{99}{112.42}{31.00}{100}
\emline{112.42}{31.00}{101}{127.42}{31.00}{102}
\emline{127.42}{31.00}{103}{127.42}{41.00}{104}
\emline{127.42}{31.00}{105}{132.42}{31.00}{106}
\emline{132.42}{31.00}{107}{132.42}{41.00}{108}
\emline{132.42}{36.00}{109}{127.42}{36.00}{110}
\emline{132.42}{31.00}{111}{147.42}{31.00}{112}
\emline{147.42}{31.00}{113}{147.42}{36.00}{114}
\emline{147.42}{31.00}{115}{152.42}{31.00}{116}
\emline{152.42}{31.00}{117}{152.42}{36.00}{118}
\emline{92.42}{36.00}{119}{92.42}{31.00}{120}
\emline{92.42}{31.00}{121}{97.42}{31.00}{122}
\emline{97.42}{31.00}{123}{97.42}{36.00}{124}
\emline{97.42}{31.00}{125}{97.42}{26.00}{126}
\emline{97.42}{26.00}{127}{112.42}{26.00}{128}
\emline{112.42}{26.00}{129}{112.42}{31.00}{130}
\emline{112.42}{26.00}{131}{122.42}{26.00}{132}
\emline{122.42}{26.00}{133}{122.42}{31.00}{134}
\emline{117.42}{31.00}{135}{117.42}{26.00}{136}
\emline{122.42}{26.00}{137}{122.42}{21.00}{138}
\emline{122.42}{21.00}{139}{137.42}{21.00}{140}
\emline{137.42}{21.00}{141}{137.42}{31.00}{142}
\emline{137.42}{21.00}{143}{142.42}{21.00}{144}
\emline{142.42}{21.00}{145}{142.42}{31.00}{146}
\emline{137.42}{26.00}{147}{142.42}{26.00}{148}
\emline{142.42}{21.00}{149}{152.42}{21.00}{150}
\emline{92.42}{21.00}{151}{97.42}{21.00}{152}
\emline{97.42}{21.00}{153}{97.42}{26.00}{154}
\emline{97.42}{21.00}{155}{107.42}{21.00}{156}
\emline{107.42}{21.00}{157}{107.42}{26.00}{158}
\emline{102.42}{26.00}{159}{102.42}{21.00}{160}
\emline{107.42}{21.00}{161}{107.42}{16.00}{162}
\emline{107.42}{16.00}{163}{132.42}{16.00}{164}
\emline{132.42}{16.00}{165}{132.42}{21.00}{166}
\emline{127.42}{21.00}{167}{127.42}{16.00}{168}
\emline{122.42}{16.00}{169}{122.42}{21.00}{170}
\emline{132.42}{16.00}{171}{132.42}{11.00}{172}
\emline{132.42}{11.00}{173}{147.42}{11.00}{174}
\emline{147.42}{11.00}{175}{147.42}{21.00}{176}
\emline{147.42}{11.00}{177}{152.42}{11.00}{178}
\emline{152.42}{11.00}{179}{152.42}{21.00}{180}
\emline{152.42}{16.00}{181}{147.42}{16.00}{182}
\emline{92.42}{21.00}{183}{92.42}{11.00}{184}
\emline{92.42}{11.00}{185}{107.42}{11.00}{186}
\emline{107.42}{11.00}{187}{107.42}{16.00}{188}
\emline{107.42}{11.00}{189}{117.42}{11.00}{190}
\emline{112.42}{16.00}{191}{112.42}{11.00}{192}
\emline{117.42}{16.00}{193}{117.42}{6.00}{194}
\emline{117.42}{6.00}{195}{142.42}{6.00}{196}
\emline{142.42}{6.00}{197}{142.42}{11.00}{198}
\emline{137.42}{11.00}{199}{137.42}{6.00}{200}
\emline{142.42}{6.00}{201}{142.42}{1.00}{202}
\emline{142.42}{1.00}{203}{152.42}{1.00}{204}
\emline{92.42}{1.00}{205}{97.42}{1.00}{206}
\emline{97.42}{1.00}{207}{97.42}{11.00}{208}
\emline{97.42}{1.00}{209}{102.42}{1.00}{210}
\emline{102.42}{1.00}{211}{102.42}{11.00}{212}
\emline{97.42}{6.00}{213}{102.42}{6.00}{214}
\emline{102.42}{1.00}{215}{117.42}{1.00}{216}
\emline{117.42}{1.00}{217}{117.42}{6.00}{218}
\emline{117.42}{1.00}{219}{127.42}{1.00}{220}
\emline{127.42}{1.00}{221}{127.42}{6.00}{222}
\emline{122.42}{6.00}{223}{122.42}{1.00}{224}
\put(97.42,56.00){\circle*{2.00}}
\put(102.42,56.00){\circle*{2.00}}
\put(117.42,56.00){\circle*{2.00}}
\put(122.42,56.00){\circle*{2.00}}
\put(107.42,46.00){\circle*{2.00}}
\put(112.42,46.00){\circle*{2.00}}
\put(127.42,46.00){\circle*{2.00}}
\put(132.42,46.00){\circle*{2.00}}
\put(117.42,36.00){\circle*{2.00}}
\put(122.42,36.00){\circle*{2.00}}
\put(137.42,36.00){\circle*{2.00}}
\put(142.42,36.00){\circle*{2.00}}
\put(127.42,26.00){\circle*{2.00}}
\put(132.42,26.00){\circle*{2.00}}
\put(147.42,26.00){\circle*{2.00}}
\put(152.42,26.00){\circle*{2.00}}
\put(137.42,16.00){\circle*{2.00}}
\put(142.42,16.00){\circle*{2.00}}
\put(147.42,6.00){\circle*{2.00}}
\put(152.42,6.00){\circle*{2.00}}
\put(92.42,41.00){\circle*{2.00}}
\put(97.42,41.00){\circle*{2.00}}
\put(102.42,31.00){\circle*{2.00}}
\put(107.42,31.00){\circle*{2.00}}
\put(112.42,21.00){\circle*{2.00}}
\put(117.42,21.00){\circle*{2.00}}
\put(122.42,11.00){\circle*{2.00}}
\put(127.42,11.00){\circle*{2.00}}
\put(132.42,1.00){\circle*{2.00}}
\put(137.42,1.00){\circle*{2.00}}
\put(132.42,61.00){\circle*{2.00}}
\put(137.42,61.00){\circle*{2.00}}
\put(92.42,26.00){\circle*{2.00}}
\put(97.42,16.00){\circle*{2.00}}
\put(102.42,16.00){\circle*{2.00}}
\put(107.42,6.00){\circle*{2.00}}
\put(112.42,6.00){\circle*{2.00}}
\put(92.42,6.00){\circle*{2.00}}
\put(142.42,51.00){\circle*{2.00}}
\put(147.42,51.00){\circle*{2.00}}
\put(152.42,41.00){\circle*{2.00}}
\emline{132.42}{6.00}{225}{132.42}{11.00}{226}
\emline{67.42}{61.00}{227}{67.42}{62.00}{228}
\emline{72.42}{61.00}{229}{72.42}{62.00}{230}
\emline{62.42}{56.00}{231}{61.42}{56.00}{232}
\emline{62.42}{51.00}{233}{61.42}{51.00}{234}
\emline{62.42}{41.00}{235}{61.42}{41.00}{236}
\emline{77.42}{41.00}{237}{77.42}{40.00}{238}
\emline{62.42}{41.00}{239}{62.42}{40.00}{240}
\emline{82.42}{41.00}{241}{82.42}{40.00}{242}
\emline{82.42}{41.00}{243}{83.42}{41.00}{244}
\emline{82.42}{51.00}{245}{83.42}{51.00}{246}
\emline{82.42}{61.00}{247}{83.42}{61.00}{248}
\emline{61.42}{61.00}{249}{62.42}{61.00}{250}
\emline{97.42}{61.00}{251}{97.42}{62.00}{252}
\emline{102.42}{61.00}{253}{102.42}{62.00}{254}
\emline{117.42}{61.00}{255}{117.42}{62.00}{256}
\emline{122.42}{61.00}{257}{122.42}{62.00}{258}
\emline{107.42}{1.00}{259}{107.42}{0.00}{260}
\emline{112.42}{1.00}{261}{112.42}{0.00}{262}
\emline{147.42}{1.00}{263}{147.42}{0.00}{264}
\emline{152.42}{1.00}{265}{152.42}{0.00}{266}
\emline{152.42}{1.00}{267}{153.42}{1.00}{268}
\emline{92.42}{1.00}{269}{91.42}{1.00}{270}
\emline{92.42}{1.00}{271}{92.42}{0.00}{272}
\emline{92.42}{61.00}{273}{91.42}{61.00}{274}
\emline{92.42}{56.00}{275}{91.42}{56.00}{276}
\emline{92.42}{46.00}{277}{91.42}{46.00}{278}
\emline{92.42}{36.00}{279}{91.42}{36.00}{280}
\emline{92.42}{31.00}{281}{91.42}{31.00}{282}
\emline{92.42}{21.00}{283}{91.42}{21.00}{284}
\emline{92.42}{16.00}{285}{91.42}{16.00}{286}
\emline{92.42}{11.00}{287}{91.42}{11.00}{288}
\emline{152.42}{11.00}{289}{153.42}{11.00}{290}
\emline{152.42}{21.00}{291}{153.42}{21.00}{292}
\emline{152.42}{31.00}{293}{153.42}{31.00}{294}
\emline{152.42}{36.00}{295}{153.42}{36.00}{296}
\emline{152.42}{46.00}{297}{153.42}{46.00}{298}
\emline{152.42}{51.00}{299}{153.42}{51.00}{300}
\emline{153.42}{61.00}{301}{153.42}{61.00}{302}
\emline{152.42}{61.00}{303}{153.42}{61.00}{304}
\emline{142.42}{61.00}{305}{142.42}{62.00}{306}
\emline{127.42}{61.00}{307}{127.42}{62.00}{308}
\end{picture}
\caption{First two even-weight cases: $M(0)$ and $M(011)$}
\end{figure}

\begin{cor}
A PTPC in a Cartesian product $C_m\times C_n$ of two cycles $C_m$ and $C_n$
exists if and only if $m$ and $n$ are
multiples of 4.
\end{cor}

\proof Let $B$ be a fixed period $B$ of a periodic doubly infinite binary sequence.
By taking the union of several contiguous copies
of $N(B)$ conforming a square or rectangle $N'$ in the Euclidean plane, or just
$N'=N(B)$, and
identifying the top and bottom sides of $N'$ as well as the left and right sides of $N'$,
it is seen that a toroidal graph of the form $C_m\times C_n$ is obtained that contains a
PTPC. Because of Theorems 1 and 4, both $m$ and $n$ must be multiples of 4. Moreover,
no other sources of PTPCs in Cartesian product of cycles exist.
\qfd

Each Cartesian product $C_m\times C_n$ of cycles $C_m$ and $C_n$ of lengths
$m$ and $n$ larger than 2, respectively,
is the target of a {\it canonical projection} graph map
$\rho$ from ${\Lambda}$ onto $C_m\times C_n$ such that
$\rho(i,j)=(i\,\,mod\,\, m,j\,\,mod\,\, n)$,
where $i\,\,mod\,\,m$ and $j\,\,mod\,\,n$ are the
remainders of dividing $i$ and $j$ respectively by $m$ and $n$. Clearly,
$\rho$ is surjective graph homomorphism.

\begin{figure}
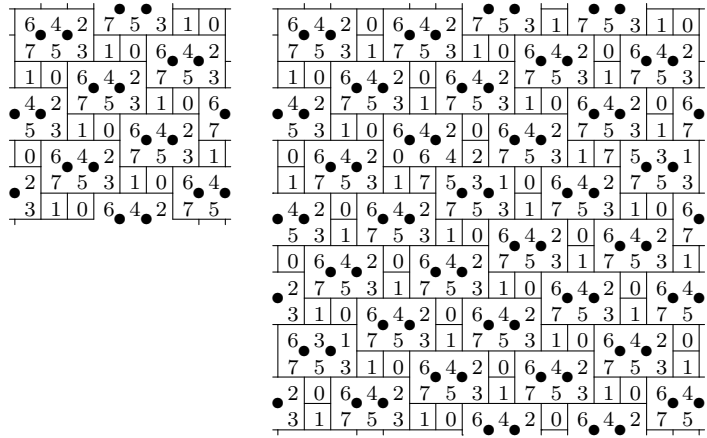

\unitlength=0.70mm
\special{em:linewidth 0.4pt}
\linethickness{0.4pt}

\caption{First two odd-weight cases: $M(1)$ and $M(01)$}
\end{figure}

\begin{cor}
A TPC in a Cartesian product $C_m\times C_n$ producing only ladder pairs of the form 21
in the inverse-image TPC in ${\Lambda}$ via $\rho$
exists if and only if $m$ and $n$ are multiples of 4.
\end{cor}

\proof
The statement of the corollary arises by combining Corollary 3 and Theorem 4.
\qfd

\section{$d$-Perfect code partitions of ${\Lambda}$ and $C_m\times C_n$}

In contrast with the uncountability of Theorem 1, there is just
one 1-perfect code $S_0$ in ${\Lambda}$, up to symmetry, which can
be taken as the sublattice of ${\Lambda}$ generated by
$\{(1,2),(2,-1)\}$. (A generalization of this fact is given in
Section 8 for integer lattices in $\R^r$, $r\ge 2$). Figure 6
partially depicts $S_0$ and its complementary graph in
${\Lambda}$. Moreover, $S_0$ does not restrict to a 1-perfect code
in any rectangular grid graph ${\Gamma}_{m,n}$ with $m$ or $n$
larger than $4$. The center $4\times 4$ grid graph
${\Gamma}_{4,4}$ in the interior of the dotted line in Figure 6
depicts the only existing 1-perfect code in a ${\Gamma}_{m,n}$
with $min\{m,n\}>2$, up to symmetry, \cite{LS}. (It is not
difficult to visualize 1-perfect codes in ${\Gamma}_{m,n}$ if
$m=1$ or 2, in the latter case with $n$ odd).

$V({\Lambda})$ admits a partition into five copies of $S_0$. Now,
the total number of perfect codes isomorphic to $S_0$ in ${\Lambda}$
is ten, of which five are tilted as in Figure 6 and five are tilted
in the other way, composing two {\it enantiomorphic} presentations
of $S_0$, (mirror images of each other). Thus, there are two
possible partitions of ${\Lambda}$ into copies of $S_0$, or of its
enantiomorphic code in ${\Lambda}$.

\begin{figure}
\unitlength=0.30mm
\special{em:linewidth 0.4pt}
\linethickness{0.4pt}
\begin{picture}(355.22,91.95)
\put(294.22,50.95){\circle*{2.00}}
\put(314.22,60.95){\circle*{2.00}}
\put(324.22,40.95){\circle*{2.00}}
\put(304.22,30.95){\circle*{2.00}}
\put(334.22,70.95){\circle*{2.00}}
\put(354.22,80.95){\circle*{2.00}}
\put(274.22,40.95){\circle*{2.00}}
\put(284.22,20.95){\circle*{2.00}}
\put(294.22,0.95){\circle*{2.00}}
\put(314.22,10.95){\circle*{2.00}}
\put(334.22,20.95){\circle*{2.00}}
\put(354.22,30.95){\circle*{2.00}}
\put(344.22,50.95){\circle*{2.00}}
\put(324.22,90.95){\circle*{2.00}}
\put(304.22,80.95){\circle*{2.00}}
\put(284.22,70.95){\circle*{2.00}}
\put(264.22,60.95){\circle*{2.00}}
\put(274.22,90.95){\circle*{2.00}}
\put(264.22,10.95){\circle*{2.00}}
\put(344.22,0.95){\circle*{2.00}}
\emline{286.22}{70.95}{1}{287.22}{70.95}{2}
\emline{288.22}{70.95}{3}{289.22}{70.95}{4}
\emline{290.22}{70.95}{5}{291.22}{70.95}{6}
\emline{292.22}{70.95}{7}{293.22}{70.95}{8}
\emline{326.22}{70.95}{9}{327.22}{70.95}{10}
\emline{328.22}{70.95}{11}{329.22}{70.95}{12}
\emline{330.22}{70.95}{13}{331.22}{70.95}{14}
\emline{332.22}{70.95}{15}{333.22}{70.95}{16}
\emline{286.22}{20.95}{17}{287.22}{20.95}{18}
\emline{288.22}{20.95}{19}{289.22}{20.95}{20}
\emline{290.22}{20.95}{21}{291.22}{20.95}{22}
\emline{292.22}{20.95}{23}{293.22}{20.95}{24}
\emline{326.22}{20.95}{25}{327.22}{20.95}{26}
\emline{328.22}{20.95}{27}{329.22}{20.95}{28}
\emline{330.22}{20.95}{29}{331.22}{20.95}{30}
\emline{332.22}{20.95}{31}{333.22}{20.95}{32}
\emline{284.22}{68.95}{33}{284.22}{67.95}{34}
\emline{284.22}{66.95}{35}{284.22}{65.95}{36}
\emline{284.22}{64.95}{37}{284.22}{63.95}{38}
\emline{284.22}{62.95}{39}{284.22}{61.95}{40}
\emline{284.22}{28.95}{41}{284.22}{27.95}{42}
\emline{284.22}{26.95}{43}{284.22}{25.95}{44}
\emline{284.22}{24.95}{45}{284.22}{23.95}{46}
\emline{284.22}{22.95}{47}{284.22}{21.95}{48}
\emline{334.22}{68.95}{49}{334.22}{67.95}{50}
\emline{334.22}{66.95}{51}{334.22}{65.95}{52}
\emline{334.22}{64.95}{53}{334.22}{63.95}{54}
\emline{334.22}{62.95}{55}{334.22}{61.95}{56}
\emline{334.22}{28.95}{57}{334.22}{27.95}{58}
\emline{334.22}{26.95}{59}{334.22}{25.95}{60}
\emline{334.22}{24.95}{61}{334.22}{23.95}{62}
\emline{334.22}{22.95}{63}{334.22}{21.95}{64}
\emline{284.22}{40.95}{65}{314.22}{40.95}{66}
\emline{304.22}{40.95}{67}{304.22}{70.95}{68}
\emline{274.22}{60.95}{69}{304.22}{60.95}{70}
\emline{264.22}{50.95}{71}{284.22}{50.95}{72}
\emline{274.22}{80.95}{73}{274.22}{50.95}{74}
\emline{264.22}{80.95}{75}{294.22}{80.95}{76}
\emline{264.22}{90.95}{77}{264.22}{70.95}{78}
\emline{264.22}{50.95}{79}{264.22}{20.95}{80}
\emline{264.22}{30.95}{81}{294.22}{30.95}{82}
\emline{274.22}{30.95}{83}{274.22}{0.95}{84}
\emline{264.22}{0.95}{85}{284.22}{0.95}{86}
\emline{284.22}{0.95}{87}{284.22}{10.95}{88}
\emline{274.22}{10.95}{89}{304.22}{10.95}{90}
\emline{304.22}{0.95}{91}{334.22}{0.95}{92}
\emline{334.22}{0.95}{93}{334.22}{10.95}{94}
\emline{324.22}{10.95}{95}{354.22}{10.95}{96}
\emline{354.22}{0.95}{97}{354.22}{20.95}{98}
\emline{344.22}{20.95}{99}{354.22}{20.95}{100}
\emline{344.22}{10.95}{101}{344.22}{40.95}{102}
\emline{354.22}{40.95}{103}{334.22}{40.95}{104}
\emline{354.22}{40.95}{105}{354.22}{70.95}{106}
\emline{344.22}{60.95}{107}{344.22}{90.95}{108}
\emline{354.22}{70.95}{109}{344.22}{70.95}{110}
\emline{344.22}{80.95}{111}{314.22}{80.95}{112}
\emline{334.22}{90.95}{113}{354.22}{90.95}{114}
\emline{314.22}{90.95}{115}{284.22}{90.95}{116}
\emline{294.22}{90.95}{117}{294.22}{60.95}{118}
\emline{284.22}{90.95}{119}{284.22}{80.95}{120}
\emline{264.22}{70.95}{121}{274.22}{70.95}{122}
\emline{294.22}{40.95}{123}{294.22}{10.95}{124}
\emline{304.22}{0.95}{125}{304.22}{20.95}{126}
\emline{314.22}{20.95}{127}{314.22}{50.95}{128}
\emline{304.22}{50.95}{129}{334.22}{50.95}{130}
\emline{324.22}{80.95}{131}{324.22}{50.95}{132}
\emline{324.22}{30.95}{133}{324.22}{0.95}{134}
\emline{344.22}{30.95}{135}{314.22}{30.95}{136}
\emline{294.22}{20.95}{137}{324.22}{20.95}{138}
\emline{334.22}{30.95}{139}{334.22}{60.95}{140}
\emline{324.22}{70.95}{141}{294.22}{70.95}{142}
\emline{284.22}{60.95}{143}{284.22}{30.95}{144}
\emline{324.22}{60.95}{145}{354.22}{60.95}{146}
\emline{314.22}{90.95}{147}{314.22}{70.95}{148}
\emline{334.22}{90.95}{149}{334.22}{80.95}{150}
\emline{264.22}{20.95}{151}{274.22}{20.95}{152}
\put(144.22,0.95){\circle{2.00}}
\put(154.22,0.95){\circle{2.00}}
\put(164.22,0.95){\circle{2.00}}
\put(174.22,0.95){\circle{2.00}}
\put(184.22,0.95){\circle{2.00}}
\put(194.22,0.95){\circle{2.00}}
\put(204.22,0.95){\circle{2.00}}
\put(214.22,0.95){\circle{2.00}}
\put(144.22,10.95){\circle{2.00}}
\put(154.22,10.95){\circle{2.00}}
\put(164.22,10.95){\circle{2.00}}
\put(174.22,10.95){\circle{2.00}}
\put(184.22,10.95){\circle{2.00}}
\put(194.22,10.95){\circle{2.00}}
\put(204.22,10.95){\circle{2.00}}
\put(214.22,10.95){\circle{2.00}}
\put(144.22,20.95){\circle{2.00}}
\put(154.22,20.95){\circle{2.00}}
\put(164.22,20.95){\circle{2.00}}
\put(174.22,20.95){\circle{2.00}}
\put(184.22,20.95){\circle{2.00}}
\put(194.22,20.95){\circle{2.00}}
\put(204.22,20.95){\circle{2.00}}
\put(214.22,20.95){\circle{2.00}}
\put(144.22,30.95){\circle{2.00}}
\put(154.22,30.95){\circle{2.00}}
\put(164.22,30.95){\circle{2.00}}
\put(174.22,30.95){\circle{2.00}}
\put(184.22,30.95){\circle{2.00}}
\put(194.22,30.95){\circle{2.00}}
\put(204.22,30.95){\circle{2.00}}
\put(214.22,30.95){\circle{2.00}}
\put(144.22,40.95){\circle{2.00}}
\put(154.22,40.95){\circle{2.00}}
\put(164.22,40.95){\circle{2.00}}
\put(174.22,40.95){\circle{2.00}}
\put(184.22,40.95){\circle{2.00}}
\put(194.22,40.95){\circle{2.00}}
\put(204.22,40.95){\circle{2.00}}
\put(214.22,40.95){\circle{2.00}}
\put(144.22,50.95){\circle{2.00}}
\put(154.22,50.95){\circle{2.00}}
\put(164.22,50.95){\circle{2.00}}
\put(174.22,50.95){\circle{2.00}}
\put(184.22,50.95){\circle{2.00}}
\put(194.22,50.95){\circle{2.00}}
\put(204.22,50.95){\circle{2.00}}
\put(214.22,50.95){\circle{2.00}}
\put(144.22,60.95){\circle{2.00}}
\put(154.22,60.95){\circle{2.00}}
\put(164.22,60.95){\circle{2.00}}
\put(174.22,60.95){\circle{2.00}}
\put(184.22,60.95){\circle{2.00}}
\put(194.22,60.95){\circle{2.00}}
\put(204.22,60.95){\circle{2.00}}
\put(214.22,60.95){\circle{2.00}}
\put(144.22,70.95){\circle{2.00}}
\put(154.22,70.95){\circle{2.00}}
\put(164.22,70.95){\circle{2.00}}
\put(174.22,70.95){\circle{2.00}}
\put(184.22,70.95){\circle{2.00}}
\put(194.22,70.95){\circle{2.00}}
\put(204.22,70.95){\circle{2.00}}
\put(214.22,70.95){\circle{2.00}}
\put(224.22,0.95){\circle{2.00}}
\put(234.22,0.95){\circle{2.00}}
\put(224.22,10.95){\circle{2.00}}
\put(234.22,10.95){\circle{2.00}}
\put(224.22,20.95){\circle{2.00}}
\put(234.22,20.95){\circle{2.00}}
\put(224.22,30.95){\circle{2.00}}
\put(234.22,30.95){\circle{2.00}}
\put(224.22,40.95){\circle{2.00}}
\put(234.22,40.95){\circle{2.00}}
\put(224.22,50.95){\circle{2.00}}
\put(234.22,50.95){\circle{2.00}}
\put(224.22,60.95){\circle{2.00}}
\put(234.22,60.95){\circle{2.00}}
\put(224.22,70.95){\circle{2.00}}
\put(234.22,70.95){\circle{2.00}}
\put(144.22,80.95){\circle{2.00}}
\put(154.22,80.95){\circle{2.00}}
\put(164.22,80.95){\circle{2.00}}
\put(174.22,80.95){\circle{2.00}}
\put(184.22,80.95){\circle{2.00}}
\put(194.22,80.95){\circle{2.00}}
\put(204.22,80.95){\circle{2.00}}
\put(214.22,80.95){\circle{2.00}}
\put(144.22,90.95){\circle{2.00}}
\put(154.22,90.95){\circle{2.00}}
\put(164.22,90.95){\circle{2.00}}
\put(174.22,90.95){\circle{2.00}}
\put(184.22,90.95){\circle{2.00}}
\put(194.22,90.95){\circle{2.00}}
\put(204.22,90.95){\circle{2.00}}
\put(214.22,90.95){\circle{2.00}}
\put(224.22,80.95){\circle{2.00}}
\put(234.22,80.95){\circle{2.00}}
\put(224.22,90.95){\circle{2.00}}
\put(234.22,90.95){\circle{2.00}}
\put(174.22,50.95){\circle*{2.00}}
\put(194.22,60.95){\circle*{2.00}}
\put(204.22,40.95){\circle*{2.00}}
\put(184.22,30.95){\circle*{2.00}}
\put(214.22,70.95){\circle*{2.00}}
\put(234.22,80.95){\circle*{2.00}}
\put(154.22,40.95){\circle*{2.00}}
\put(164.22,20.95){\circle*{2.00}}
\put(174.22,0.95){\circle*{2.00}}
\put(194.22,10.95){\circle*{2.00}}
\put(214.22,20.95){\circle*{2.00}}
\put(234.22,30.95){\circle*{2.00}}
\put(224.22,50.95){\circle*{2.00}}
\put(204.22,90.95){\circle*{2.00}}
\put(184.22,80.95){\circle*{2.00}}
\put(164.22,70.95){\circle*{2.00}}
\put(144.22,60.95){\circle*{2.00}}
\put(154.22,90.95){\circle*{2.00}}
\put(144.22,10.95){\circle*{2.00}}
\put(224.22,0.95){\circle*{2.00}}
\emline{166.22}{70.95}{153}{167.22}{70.95}{154}
\emline{168.22}{70.95}{155}{169.22}{70.95}{156}
\emline{170.22}{70.95}{157}{171.22}{70.95}{158}
\emline{172.22}{70.95}{159}{173.22}{70.95}{160}
\emline{176.22}{70.95}{161}{177.22}{70.95}{162}
\emline{178.22}{70.95}{163}{179.22}{70.95}{164}
\emline{180.22}{70.95}{165}{181.22}{70.95}{166}
\emline{182.22}{70.95}{167}{183.22}{70.95}{168}
\emline{186.22}{70.95}{169}{187.22}{70.95}{170}
\emline{188.22}{70.95}{171}{189.22}{70.95}{172}
\emline{190.22}{70.95}{173}{191.22}{70.95}{174}
\emline{192.22}{70.95}{175}{193.22}{70.95}{176}
\emline{196.22}{70.95}{177}{197.22}{70.95}{178}
\emline{198.22}{70.95}{179}{199.22}{70.95}{180}
\emline{200.22}{70.95}{181}{201.22}{70.95}{182}
\emline{202.22}{70.95}{183}{203.22}{70.95}{184}
\emline{206.22}{70.95}{185}{207.22}{70.95}{186}
\emline{208.22}{70.95}{187}{209.22}{70.95}{188}
\emline{210.22}{70.95}{189}{211.22}{70.95}{190}
\emline{212.22}{70.95}{191}{213.22}{70.95}{192}
\emline{166.22}{20.95}{193}{167.22}{20.95}{194}
\emline{168.22}{20.95}{195}{169.22}{20.95}{196}
\emline{170.22}{20.95}{197}{171.22}{20.95}{198}
\emline{172.22}{20.95}{199}{173.22}{20.95}{200}
\emline{176.22}{20.95}{201}{177.22}{20.95}{202}
\emline{178.22}{20.95}{203}{179.22}{20.95}{204}
\emline{180.22}{20.95}{205}{181.22}{20.95}{206}
\emline{182.22}{20.95}{207}{183.22}{20.95}{208}
\emline{186.22}{20.95}{209}{187.22}{20.95}{210}
\emline{188.22}{20.95}{211}{189.22}{20.95}{212}
\emline{190.22}{20.95}{213}{191.22}{20.95}{214}
\emline{192.22}{20.95}{215}{193.22}{20.95}{216}
\emline{196.22}{20.95}{217}{197.22}{20.95}{218}
\emline{198.22}{20.95}{219}{199.22}{20.95}{220}
\emline{200.22}{20.95}{221}{201.22}{20.95}{222}
\emline{202.22}{20.95}{223}{203.22}{20.95}{224}
\emline{206.22}{20.95}{225}{207.22}{20.95}{226}
\emline{208.22}{20.95}{227}{209.22}{20.95}{228}
\emline{210.22}{20.95}{229}{211.22}{20.95}{230}
\emline{212.22}{20.95}{231}{213.22}{20.95}{232}
\emline{164.22}{68.95}{233}{164.22}{67.95}{234}
\emline{164.22}{66.95}{235}{164.22}{65.95}{236}
\emline{164.22}{64.95}{237}{164.22}{63.95}{238}
\emline{164.22}{62.95}{239}{164.22}{61.95}{240}
\emline{164.22}{58.95}{241}{164.22}{57.95}{242}
\emline{164.22}{56.95}{243}{164.22}{55.95}{244}
\emline{164.22}{54.95}{245}{164.22}{53.95}{246}
\emline{164.22}{52.95}{247}{164.22}{51.95}{248}
\emline{164.22}{48.95}{249}{164.22}{47.95}{250}
\emline{164.22}{46.95}{251}{164.22}{45.95}{252}
\emline{164.22}{44.95}{253}{164.22}{43.95}{254}
\emline{164.22}{42.95}{255}{164.22}{41.95}{256}
\emline{164.22}{38.95}{257}{164.22}{37.95}{258}
\emline{164.22}{36.95}{259}{164.22}{35.95}{260}
\emline{164.22}{34.95}{261}{164.22}{33.95}{262}
\emline{164.22}{32.95}{263}{164.22}{31.95}{264}
\emline{164.22}{28.95}{265}{164.22}{27.95}{266}
\emline{164.22}{26.95}{267}{164.22}{25.95}{268}
\emline{164.22}{24.95}{269}{164.22}{23.95}{270}
\emline{164.22}{22.95}{271}{164.22}{21.95}{272}
\emline{214.22}{68.95}{273}{214.22}{67.95}{274}
\emline{214.22}{66.95}{275}{214.22}{65.95}{276}
\emline{214.22}{64.95}{277}{214.22}{63.95}{278}
\emline{214.22}{62.95}{279}{214.22}{61.95}{280}
\emline{214.22}{58.95}{281}{214.22}{57.95}{282}
\emline{214.22}{56.95}{283}{214.22}{55.95}{284}
\emline{214.22}{54.95}{285}{214.22}{53.95}{286}
\emline{214.22}{52.95}{287}{214.22}{51.95}{288}
\emline{214.22}{48.95}{289}{214.22}{47.95}{290}
\emline{214.22}{46.95}{291}{214.22}{45.95}{292}
\emline{214.22}{44.95}{293}{214.22}{43.95}{294}
\emline{214.22}{42.95}{295}{214.22}{41.95}{296}
\emline{214.22}{38.95}{297}{214.22}{37.95}{298}
\emline{214.22}{36.95}{299}{214.22}{35.95}{300}
\emline{214.22}{34.95}{301}{214.22}{33.95}{302}
\emline{214.22}{32.95}{303}{214.22}{31.95}{304}
\emline{214.22}{28.95}{305}{214.22}{27.95}{306}
\emline{214.22}{26.95}{307}{214.22}{25.95}{308}
\emline{214.22}{24.95}{309}{214.22}{23.95}{310}
\emline{214.22}{22.95}{311}{214.22}{21.95}{312}
\emline{144.22}{10.95}{313}{164.22}{20.95}{314}
\emline{164.22}{20.95}{315}{154.22}{40.95}{316}
\emline{154.22}{40.95}{317}{144.22}{60.95}{318}
\emline{144.22}{60.95}{319}{164.22}{70.95}{320}
\emline{164.22}{70.95}{321}{154.22}{90.95}{322}
\emline{154.22}{40.95}{323}{174.22}{50.95}{324}
\emline{174.22}{50.95}{325}{164.22}{70.95}{326}
\emline{164.22}{70.95}{327}{184.22}{80.95}{328}
\emline{184.22}{80.95}{329}{204.22}{90.95}{330}
\emline{204.22}{90.95}{331}{214.22}{70.95}{332}
\emline{214.22}{70.95}{333}{194.22}{60.95}{334}
\emline{194.22}{60.95}{335}{204.22}{40.95}{336}
\emline{204.22}{40.95}{337}{184.22}{30.95}{338}
\emline{184.22}{30.95}{339}{194.22}{10.95}{340}
\emline{194.22}{10.95}{341}{214.22}{20.95}{342}
\emline{214.22}{20.95}{343}{234.22}{30.95}{344}
\emline{234.22}{30.95}{345}{224.22}{50.95}{346}
\emline{224.22}{50.95}{347}{214.22}{70.95}{348}
\emline{214.22}{70.95}{349}{234.22}{80.95}{350}
\emline{184.22}{80.95}{351}{194.22}{60.95}{352}
\emline{194.22}{60.95}{353}{174.22}{50.95}{354}
\emline{174.22}{50.95}{355}{184.22}{30.95}{356}
\emline{184.22}{30.95}{357}{164.22}{20.95}{358}
\emline{164.22}{20.95}{359}{174.22}{0.95}{360}
\emline{174.22}{0.95}{361}{194.22}{10.95}{362}
\emline{224.22}{50.95}{363}{204.22}{40.95}{364}
\emline{204.22}{40.95}{365}{214.22}{20.95}{366}
\emline{214.22}{20.95}{367}{224.22}{0.95}{368}
\end{picture}
\caption{1-Perfect code in ${\Lambda}$ and a cutout for $C_5\times C_5$}
\end{figure}
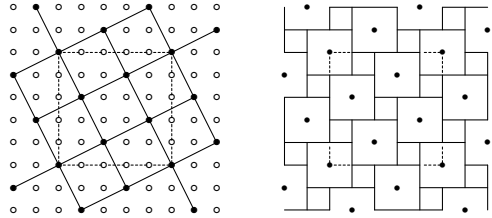

Moreover, there are 1-perfect codes in toroidal graphs,
specifically Cartesian products $C_{5k}\times C_{5\ell}$,
for $0<k,\ell\in\Z$, obtained from $S_0$ and having cardinality $5k\ell$.
(The existence of such a code would follow from Theorem 2.5 of \cite{LS},
but look at the remark after Theorem 11 in Section 6.)

An example of such 1-perfect code in $C_5\times C_5$ can
be visualized in Figure 6, where the dotted lines delineate the
boundary of a cutout of the (flat) torus involved; the left
side of the figure represents the minimum-distance graph of $S_0$,
which is the restriction of the power graph $({\Lambda})^3$ to $S_0$;
the right side represents
the complementary graph of $S_0$ in ${\Lambda}$. No other 1-perfect codes
in Cartesian products of two cycles exist.

\begin{thm}
There exists a toroidal graph $C_m\times C_n$ having a
1-perfect code partition ${\mathcal
S}_0=\{S_0^0=S_0,S_0^1,\ldots,S_0^4\}$ if and only if $m$ and $n$
are multiples of 5. Each component code $S_0^i$ of ${\mathcal
S}_0$, which is a translate of the sublattice $S_0$, has
cardinality $mn/5$ and cannot be obtained by side identifications
from 1-perfect codes in any rectangular grid graph.
\end{thm}

Theorem 7 is a particular case of the following result.
Recall from \cite{biggs} that for $d\geq 1$,
a $d$-perfect code in a regular graph $G$ is a vertex subset $S$ such that
every vertex $v$ in $G\setminus S$ is at distance $\leq d$ from just one vertex $w(v)$
of $S$. Notice however that no $G\setminus S$ is 3-regular,
for any $d$-perfect code $S$ in a 4-regular graph $G$, when $d>1$, contrary to the observation
previous to Theorem 1. Figure 7 illustrates this point.

A $d$-perfect code $S$ in
${\Lambda}$ can be found easily, following inductively a construction whose first
two steps are illustrated in Figures 6 and 7, for $d=1,2$ respectively. Given such an $S$,
the vertices of the $d$-neighborhood $N_d(v)$ of any vertex $v$ of $S$ are labeled
as follows: the top element of $N(v)$ gets label 1;
the vertices of the longest horizontal
path of $N(v)$, that is the $(d+1)$-th horizontal path,
of length $2d+1$, get the subsequent integer labels from left to right:
$2,3,\ldots, 2d+2$;
the bottom vertex of $N(v)$ gets label $2d+3$; if there are still horizontal paths
with unlabeled vertices in $N(v)$, then the $d$-th and $2d$-th
horizontal paths of $N(v)$ receive
the subsequent integer labels increasingly from left to right;
and so on, if necessary, with the $(d-1)$-th and $(2d-1)$-th horizontal paths, etc.

In other words, by denoting the rows of vertices of $N(v)$ from top to bottom
successively
as row 1, row 2, $\ldots$, row $(2d+1)$, and by considering this row
order cyclically mod $2d+1$, we label with $1,2,\ldots,2d^2+2d+1=q$ subsequently from left to
right the vertices in the rows
$1,1+d,1+2d,\ldots,1+(2d)d=d+2$ mod $2d+1$, in that order.
For Figures 6 and 7, this looks respective as follows:
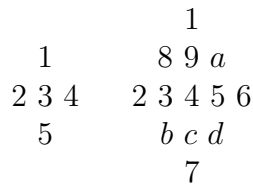
\begin{figure}
\unitlength=0.40mm
\special{em:linewidth 0.4pt}
\linethickness{0.4pt}
\begin{picture}(225.75,78.00)
\put(184.75,27.00){\circle{2.00}}
\put(184.75,32.00){\circle{2.00}}
\put(184.75,22.00){\circle{2.00}}
\put(189.75,27.00){\circle{2.00}}
\put(179.75,27.00){\circle{2.00}}
\put(184.75,27.00){\circle*{2.00}}
\emline{179.75}{37.00}{1}{189.75}{37.00}{2}
\emline{189.75}{37.00}{3}{189.75}{32.00}{4}
\emline{189.75}{32.00}{5}{194.75}{32.00}{6}
\emline{194.75}{32.00}{7}{194.75}{22.00}{8}
\emline{194.75}{22.00}{9}{189.75}{22.00}{10}
\emline{189.75}{22.00}{11}{189.75}{17.00}{12}
\emline{189.75}{17.00}{13}{179.75}{17.00}{14}
\emline{179.75}{17.00}{15}{179.75}{22.00}{16}
\emline{179.75}{22.00}{17}{174.75}{22.00}{18}
\emline{174.75}{22.00}{19}{174.75}{32.00}{20}
\emline{174.75}{32.00}{21}{179.75}{32.00}{22}
\emline{179.75}{32.00}{23}{179.75}{37.00}{24}
\emline{194.75}{22.00}{25}{199.75}{22.00}{26}
\emline{199.75}{22.00}{27}{199.75}{27.00}{28}
\emline{199.75}{27.00}{29}{194.75}{27.00}{30}
\emline{199.75}{22.00}{31}{199.75}{17.00}{32}
\emline{199.75}{17.00}{33}{204.75}{17.00}{34}
\emline{204.75}{17.00}{35}{204.75}{7.00}{36}
\emline{204.75}{7.00}{37}{199.75}{7.00}{38}
\emline{199.75}{7.00}{39}{199.75}{2.00}{40}
\emline{199.75}{2.00}{41}{189.75}{2.00}{42}
\emline{189.75}{2.00}{43}{189.75}{7.00}{44}
\emline{189.75}{7.00}{45}{184.75}{7.00}{46}
\emline{184.75}{7.00}{47}{184.75}{17.00}{48}
\emline{199.75}{27.00}{49}{204.75}{27.00}{50}
\emline{204.75}{27.00}{51}{204.75}{32.00}{52}
\emline{204.75}{32.00}{53}{209.75}{32.00}{54}
\emline{209.75}{32.00}{55}{209.75}{42.00}{56}
\emline{209.75}{42.00}{57}{204.75}{42.00}{58}
\emline{204.75}{42.00}{59}{204.75}{47.00}{60}
\emline{204.75}{47.00}{61}{194.75}{47.00}{62}
\emline{194.75}{47.00}{63}{194.75}{42.00}{64}
\emline{194.75}{42.00}{65}{189.75}{42.00}{66}
\emline{189.75}{42.00}{67}{189.75}{37.00}{68}
\emline{209.75}{32.00}{69}{214.75}{32.00}{70}
\emline{214.75}{32.00}{71}{214.75}{27.00}{72}
\emline{214.75}{27.00}{73}{219.75}{27.00}{74}
\emline{219.75}{27.00}{75}{219.75}{17.00}{76}
\emline{219.75}{17.00}{77}{214.75}{17.00}{78}
\emline{214.75}{17.00}{79}{214.75}{12.00}{80}
\emline{214.75}{12.00}{81}{204.75}{12.00}{82}
\emline{189.75}{42.00}{83}{184.75}{42.00}{84}
\emline{184.75}{42.00}{85}{184.75}{37.00}{86}
\emline{174.75}{32.00}{87}{169.75}{32.00}{88}
\emline{169.75}{32.00}{89}{169.75}{27.00}{90}
\emline{169.75}{27.00}{91}{174.75}{27.00}{92}
\emline{179.75}{17.00}{93}{179.75}{12.00}{94}
\emline{179.75}{12.00}{95}{184.75}{12.00}{96}
\emline{164.75}{27.00}{97}{164.75}{22.00}{98}
\emline{164.75}{22.00}{99}{159.75}{22.00}{100}
\emline{159.75}{22.00}{101}{159.75}{12.00}{102}
\emline{159.75}{12.00}{103}{164.75}{12.00}{104}
\emline{164.75}{12.00}{105}{164.75}{7.00}{106}
\emline{164.75}{7.00}{107}{174.75}{7.00}{108}
\emline{174.75}{7.00}{109}{174.75}{12.00}{110}
\emline{174.75}{12.00}{111}{179.75}{12.00}{112}
\emline{164.75}{27.00}{113}{169.75}{27.00}{114}
\emline{169.75}{32.00}{115}{169.75}{37.00}{116}
\emline{169.75}{37.00}{117}{164.75}{37.00}{118}
\emline{164.75}{37.00}{119}{164.75}{47.00}{120}
\emline{164.75}{47.00}{121}{169.75}{47.00}{122}
\emline{169.75}{47.00}{123}{169.75}{52.00}{124}
\emline{169.75}{52.00}{125}{179.75}{52.00}{126}
\emline{179.75}{52.00}{127}{179.75}{47.00}{128}
\emline{179.75}{47.00}{129}{184.75}{47.00}{130}
\emline{184.75}{47.00}{131}{184.75}{42.00}{132}
\emline{199.75}{52.00}{133}{204.75}{52.00}{134}
\emline{204.75}{52.00}{135}{204.75}{47.00}{136}
\emline{209.75}{37.00}{137}{214.75}{37.00}{138}
\emline{214.75}{37.00}{139}{214.75}{32.00}{140}
\emline{204.75}{7.00}{141}{209.75}{7.00}{142}
\emline{209.75}{7.00}{143}{209.75}{12.00}{144}
\emline{169.75}{7.00}{145}{169.75}{2.00}{146}
\emline{169.75}{2.00}{147}{164.75}{2.00}{148}
\emline{164.75}{2.00}{149}{164.75}{7.00}{150}
\emline{209.75}{7.00}{151}{209.75}{2.00}{152}
\emline{209.75}{2.00}{153}{214.75}{2.00}{154}
\emline{164.75}{42.00}{155}{159.75}{42.00}{156}
\emline{159.75}{42.00}{157}{159.75}{47.00}{158}
\emline{159.75}{47.00}{159}{164.75}{47.00}{160}
\put(209.75,22.00){\circle{2.00}}
\put(209.75,27.00){\circle{2.00}}
\put(209.75,17.00){\circle{2.00}}
\put(214.75,22.00){\circle{2.00}}
\put(204.75,22.00){\circle{2.00}}
\put(209.75,22.00){\circle*{2.00}}
\put(169.75,17.00){\circle{2.00}}
\put(169.75,22.00){\circle{2.00}}
\put(169.75,12.00){\circle{2.00}}
\put(174.75,17.00){\circle{2.00}}
\put(164.75,17.00){\circle{2.00}}
\put(169.75,17.00){\circle*{2.00}}
\put(194.75,12.00){\circle{2.00}}
\put(194.75,17.00){\circle{2.00}}
\put(194.75,7.00){\circle{2.00}}
\put(199.75,12.00){\circle{2.00}}
\put(189.75,12.00){\circle{2.00}}
\put(194.75,12.00){\circle*{2.00}}
\put(174.75,42.00){\circle{2.00}}
\put(174.75,47.00){\circle{2.00}}
\put(174.75,37.00){\circle{2.00}}
\put(179.75,42.00){\circle{2.00}}
\put(169.75,42.00){\circle{2.00}}
\put(174.75,42.00){\circle*{2.00}}
\put(199.75,37.00){\circle{2.00}}
\put(199.75,42.00){\circle{2.00}}
\put(199.75,32.00){\circle{2.00}}
\put(204.75,37.00){\circle{2.00}}
\put(194.75,37.00){\circle{2.00}}
\put(199.75,37.00){\circle*{2.00}}
\put(159.75,7.00){\circle{2.00}}
\emline{164.75}{2.00}{161}{159.75}{2.00}{162}
\emline{159.75}{47.00}{163}{159.75}{52.00}{164}
\emline{214.75}{37.00}{165}{219.75}{37.00}{166}
\emline{219.75}{37.00}{167}{219.75}{42.00}{168}
\emline{204.75}{52.00}{169}{209.75}{52.00}{170}
\put(209.75,47.00){\circle{2.00}}
\put(214.75,42.00){\circle{2.00}}
\put(214.75,47.00){\circle*{2.00}}
\put(214.75,52.00){\circle{2.00}}
\put(219.75,47.00){\circle{2.00}}
\put(184.75,52.00){\circle{2.00}}
\put(189.75,47.00){\circle{2.00}}
\put(189.75,52.00){\circle*{2.00}}
\put(194.75,52.00){\circle{2.00}}
\put(164.75,52.00){\circle{2.00}}
\put(179.75,2.00){\circle{2.00}}
\put(179.75,7.00){\circle{2.00}}
\put(184.75,2.00){\circle{2.00}}
\put(174.75,2.00){\circle{2.00}}
\put(179.75,2.00){\circle*{2.00}}
\put(219.75,7.00){\circle{2.00}}
\put(219.75,12.00){\circle{2.00}}
\put(219.75,2.00){\circle{2.00}}
\put(214.75,7.00){\circle{2.00}}
\put(219.75,7.00){\circle*{2.00}}
\put(204.75,2.00){\circle{2.00}}
\put(219.75,32.00){\circle{2.00}}
\put(224.75,32.00){\circle*{2.00}}
\put(224.75,37.00){\circle{2.00}}
\put(224.75,27.00){\circle{2.00}}
\emline{224.75}{42.00}{171}{224.75}{52.00}{172}
\emline{224.75}{52.00}{173}{219.75}{52.00}{174}
\emline{219.75}{52.00}{175}{219.75}{57.00}{176}
\emline{219.75}{57.00}{177}{209.75}{57.00}{178}
\emline{209.75}{57.00}{179}{209.75}{52.00}{180}
\put(189.75,57.00){\circle{2.00}}
\emline{199.75}{57.00}{181}{194.75}{57.00}{182}
\emline{194.75}{57.00}{183}{194.75}{62.00}{184}
\emline{194.75}{62.00}{185}{184.75}{62.00}{186}
\emline{184.75}{62.00}{187}{184.75}{57.00}{188}
\emline{184.75}{57.00}{189}{179.75}{57.00}{190}
\emline{179.75}{52.00}{191}{179.75}{57.00}{192}
\emline{174.75}{52.00}{193}{174.75}{57.00}{194}
\emline{174.75}{57.00}{195}{179.75}{57.00}{196}
\put(169.75,57.00){\circle{2.00}}
\put(164.75,57.00){\circle*{2.00}}
\put(164.75,62.00){\circle{2.00}}
\put(159.75,57.00){\circle{2.00}}
\emline{159.75}{52.00}{197}{154.75}{52.00}{198}
\emline{154.75}{52.00}{199}{154.75}{62.00}{200}
\emline{154.75}{62.00}{201}{159.75}{62.00}{202}
\emline{159.75}{62.00}{203}{159.75}{67.00}{204}
\emline{159.75}{67.00}{205}{169.75}{67.00}{206}
\emline{169.75}{67.00}{207}{169.75}{62.00}{208}
\emline{169.75}{62.00}{209}{174.75}{62.00}{210}
\emline{174.75}{57.00}{211}{174.75}{62.00}{212}
\put(179.75,62.00){\circle{2.00}}
\put(179.75,67.00){\circle*{2.00}}
\put(174.75,67.00){\circle{2.00}}
\put(184.75,67.00){\circle{2.00}}
\put(179.75,72.00){\circle{2.00}}
\emline{189.75}{62.00}{213}{189.75}{72.00}{214}
\emline{189.75}{72.00}{215}{184.75}{72.00}{216}
\emline{184.75}{72.00}{217}{184.75}{77.00}{218}
\emline{184.75}{77.00}{219}{174.75}{77.00}{220}
\emline{174.75}{77.00}{221}{174.75}{72.00}{222}
\emline{174.75}{72.00}{223}{169.75}{72.00}{224}
\emline{169.75}{67.00}{225}{169.75}{72.00}{226}
\emline{159.75}{42.00}{227}{154.75}{42.00}{228}
\emline{149.75}{27.00}{229}{154.75}{27.00}{230}
\emline{154.75}{27.00}{231}{154.75}{22.00}{232}
\emline{159.75}{22.00}{233}{154.75}{22.00}{234}
\put(154.75,47.00){\circle{2.00}}
\put(149.75,47.00){\circle*{2.00}}
\put(149.75,42.00){\circle{2.00}}
\put(149.75,52.00){\circle{2.00}}
\emline{219.75}{17.00}{235}{224.75}{17.00}{236}
\emline{224.75}{17.00}{237}{224.75}{22.00}{238}
\emline{224.75}{17.00}{239}{224.75}{12.00}{240}
\put(204.75,57.00){\circle{2.00}}
\put(204.75,62.00){\circle*{2.00}}
\put(199.75,62.00){\circle{2.00}}
\put(209.75,62.00){\circle{2.00}}
\put(204.75,67.00){\circle{2.00}}
\emline{214.75}{57.00}{241}{214.75}{67.00}{242}
\emline{214.75}{67.00}{243}{209.75}{67.00}{244}
\emline{209.75}{67.00}{245}{209.75}{72.00}{246}
\emline{209.75}{72.00}{247}{199.75}{72.00}{248}
\emline{199.75}{72.00}{249}{199.75}{67.00}{250}
\emline{199.75}{67.00}{251}{194.75}{67.00}{252}
\emline{194.75}{67.00}{253}{194.75}{62.00}{254}
\emline{194.75}{67.00}{255}{189.75}{67.00}{256}
\put(194.75,72.00){\circle{2.00}}
\put(194.75,77.00){\circle*{2.00}}
\put(189.75,77.00){\circle{2.00}}
\put(199.75,77.00){\circle{2.00}}
\emline{204.75}{72.00}{257}{204.75}{77.00}{258}
\emline{204.75}{77.00}{259}{209.75}{77.00}{260}
\emline{209.75}{77.00}{261}{209.75}{72.00}{262}
\emline{169.75}{72.00}{263}{164.75}{72.00}{264}
\emline{164.75}{72.00}{265}{164.75}{67.00}{266}
\put(214.75,72.00){\circle{2.00}}
\emline{214.75}{62.00}{267}{219.75}{62.00}{268}
\emline{219.75}{62.00}{269}{219.75}{57.00}{270}
\put(219.75,67.00){\circle{2.00}}
\put(219.75,72.00){\circle*{2.00}}
\put(219.75,77.00){\circle{2.00}}
\put(224.75,72.00){\circle{2.00}}
\emline{219.75}{62.00}{271}{224.75}{62.00}{272}
\emline{224.75}{62.00}{273}{224.75}{67.00}{274}
\put(224.75,57.00){\circle{2.00}}
\put(154.75,7.00){\circle*{2.00}}
\put(154.75,12.00){\circle{2.00}}
\put(154.75,2.00){\circle{2.00}}
\put(149.75,7.00){\circle{2.00}}
\emline{159.75}{17.00}{275}{154.75}{17.00}{276}
\emline{154.75}{17.00}{277}{154.75}{22.00}{278}
\emline{154.75}{17.00}{279}{149.75}{17.00}{280}
\emline{149.75}{17.00}{281}{149.75}{12.00}{282}
\put(149.75,22.00){\circle{2.00}}
\emline{154.75}{37.00}{283}{154.75}{42.00}{284}
\put(164.75,32.00){\circle{2.00}}
\put(159.75,37.00){\circle{2.00}}
\put(154.75,32.00){\circle{2.00}}
\put(159.75,27.00){\circle{2.00}}
\put(159.75,32.00){\circle*{2.00}}
\emline{149.75}{37.00}{285}{149.75}{27.00}{286}
\emline{149.75}{57.00}{287}{149.75}{62.00}{288}
\emline{149.75}{62.00}{289}{154.75}{62.00}{290}
\put(154.75,67.00){\circle{2.00}}
\put(159.75,72.00){\circle{2.00}}
\put(154.75,77.00){\circle{2.00}}
\put(149.75,72.00){\circle{2.00}}
\put(154.75,72.00){\circle*{2.00}}
\emline{164.75}{77.00}{291}{159.75}{77.00}{292}
\emline{164.75}{77.00}{293}{164.75}{72.00}{294}
\put(169.75,77.00){\circle{2.00}}
\emline{149.75}{67.00}{295}{149.75}{62.00}{296}
\put(224.75,7.00){\circle{2.00}}
\emline{149.75}{2.00}{297}{149.75}{1.00}{298}
\emline{149.75}{2.00}{299}{148.75}{2.00}{300}
\emline{159.75}{2.00}{301}{159.75}{1.00}{302}
\emline{169.75}{2.00}{303}{169.75}{1.00}{304}
\emline{189.75}{2.00}{305}{189.75}{1.00}{306}
\emline{194.75}{2.00}{307}{194.75}{1.00}{308}
\emline{214.75}{2.00}{309}{214.75}{1.00}{310}
\emline{224.75}{2.00}{311}{224.75}{1.00}{312}
\emline{224.75}{2.00}{313}{225.75}{2.00}{314}
\emline{224.75}{12.00}{315}{225.75}{12.00}{316}
\emline{224.75}{22.00}{317}{225.75}{22.00}{318}
\emline{224.75}{42.00}{319}{225.75}{42.00}{320}
\emline{224.75}{47.00}{321}{225.75}{47.00}{322}
\emline{224.75}{67.00}{323}{225.75}{67.00}{324}
\emline{149.75}{77.00}{325}{148.75}{77.00}{326}
\emline{149.75}{67.00}{327}{148.75}{67.00}{328}
\emline{149.75}{57.00}{329}{148.75}{57.00}{330}
\emline{149.75}{37.00}{331}{148.75}{37.00}{332}
\emline{149.75}{32.00}{333}{148.75}{32.00}{334}
\emline{149.75}{12.00}{335}{148.75}{12.00}{336}
\emline{179.75}{77.00}{337}{179.75}{78.00}{338}
\emline{184.75}{77.00}{339}{184.75}{78.00}{340}
\emline{204.75}{77.00}{341}{204.75}{78.00}{342}
\emline{209.75}{77.00}{343}{214.75}{77.00}{344}
\emline{214.75}{77.00}{345}{214.75}{78.00}{346}
\emline{224.75}{77.00}{347}{224.75}{78.00}{348}
\emline{224.75}{77.00}{349}{225.75}{77.00}{350}
\emline{159.75}{77.00}{351}{159.75}{78.00}{352}
\emline{149.75}{77.00}{353}{149.75}{78.00}{354}
\emline{149.75}{57.00}{355}{154.75}{57.00}{356}
\emline{149.75}{37.00}{357}{154.75}{37.00}{358}
\emline{219.75}{22.00}{359}{224.75}{22.00}{360}
\emline{224.75}{42.00}{361}{219.75}{42.00}{362}
\emline{155.75}{72.00}{363}{156.75}{72.00}{364}
\emline{157.75}{72.00}{365}{158.75}{72.00}{366}
\emline{160.75}{72.00}{367}{161.75}{72.00}{368}
\emline{162.75}{72.00}{369}{163.75}{72.00}{370}
\emline{175.75}{72.00}{371}{176.75}{72.00}{372}
\emline{177.75}{72.00}{373}{178.75}{72.00}{374}
\emline{180.75}{72.00}{375}{181.75}{72.00}{376}
\emline{182.75}{72.00}{377}{183.75}{72.00}{378}
\emline{210.75}{72.00}{379}{211.75}{72.00}{380}
\emline{212.75}{72.00}{381}{213.75}{72.00}{382}
\emline{215.75}{72.00}{383}{216.75}{72.00}{384}
\emline{217.75}{72.00}{385}{218.75}{72.00}{386}
\emline{219.75}{71.00}{387}{219.75}{70.00}{388}
\emline{219.75}{69.00}{389}{219.75}{68.00}{390}
\emline{219.75}{66.00}{391}{219.75}{65.00}{392}
\emline{219.75}{64.00}{393}{219.75}{63.00}{394}
\emline{219.75}{51.00}{395}{219.75}{50.00}{396}
\emline{219.75}{49.00}{397}{219.75}{48.00}{398}
\emline{219.75}{46.00}{399}{219.75}{45.00}{400}
\emline{219.75}{44.00}{401}{219.75}{43.00}{402}
\emline{219.75}{36.00}{403}{219.75}{35.00}{404}
\emline{219.75}{34.00}{405}{219.75}{33.00}{406}
\emline{219.75}{31.00}{407}{219.75}{30.00}{408}
\emline{219.75}{29.00}{409}{219.75}{28.00}{410}
\emline{219.75}{16.00}{411}{219.75}{15.00}{412}
\emline{219.75}{14.00}{413}{219.75}{13.00}{414}
\emline{219.75}{11.00}{415}{219.75}{10.00}{416}
\emline{219.75}{9.00}{417}{219.75}{8.00}{418}
\emline{155.75}{7.00}{419}{156.75}{7.00}{420}
\emline{157.75}{7.00}{421}{158.75}{7.00}{422}
\emline{160.75}{7.00}{423}{161.75}{7.00}{424}
\emline{162.75}{7.00}{425}{163.75}{7.00}{426}
\emline{175.75}{7.00}{427}{176.75}{7.00}{428}
\emline{177.75}{7.00}{429}{178.75}{7.00}{430}
\emline{180.75}{7.00}{431}{181.75}{7.00}{432}
\emline{182.75}{7.00}{433}{183.75}{7.00}{434}
\emline{210.75}{7.00}{435}{211.75}{7.00}{436}
\emline{212.75}{7.00}{437}{213.75}{7.00}{438}
\emline{215.75}{7.00}{439}{216.75}{7.00}{440}
\emline{217.75}{7.00}{441}{218.75}{7.00}{442}
\emline{190.75}{72.00}{443}{191.75}{72.00}{444}
\emline{192.75}{72.00}{445}{193.75}{72.00}{446}
\emline{195.75}{72.00}{447}{196.75}{72.00}{448}
\emline{197.75}{72.00}{449}{198.75}{72.00}{450}
\emline{190.75}{7.00}{451}{191.75}{7.00}{452}
\emline{192.75}{7.00}{453}{193.75}{7.00}{454}
\emline{195.75}{7.00}{455}{196.75}{7.00}{456}
\emline{197.75}{7.00}{457}{198.75}{7.00}{458}
\emline{154.75}{71.00}{459}{154.75}{70.00}{460}
\emline{154.75}{69.00}{461}{154.75}{68.00}{462}
\emline{154.75}{66.00}{463}{154.75}{65.00}{464}
\emline{154.75}{64.00}{465}{154.75}{63.00}{466}
\emline{154.75}{51.00}{467}{154.75}{50.00}{468}
\emline{154.75}{49.00}{469}{154.75}{48.00}{470}
\emline{154.75}{46.00}{471}{154.75}{45.00}{472}
\emline{154.75}{44.00}{473}{154.75}{43.00}{474}
\emline{154.75}{36.00}{475}{154.75}{35.00}{476}
\emline{154.75}{34.00}{477}{154.75}{33.00}{478}
\emline{154.75}{31.00}{479}{154.75}{30.00}{480}
\emline{154.75}{29.00}{481}{154.75}{28.00}{482}
\emline{154.75}{16.00}{483}{154.75}{15.00}{484}
\emline{154.75}{14.00}{485}{154.75}{13.00}{486}
\emline{154.75}{11.00}{487}{154.75}{10.00}{488}
\emline{154.75}{9.00}{489}{154.75}{8.00}{490}
\emline{199.75}{57.00}{491}{199.75}{47.00}{492}
\end{picture}
\caption{2-Perfect code in ${\Lambda}$ and a cutout for $C_{13}\times C_{13}$}
\end{figure}
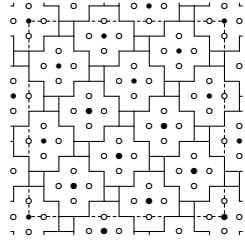
$$\begin{array}{ccc}
                   & & 1  \\
1                  & & 8\,\,9\,\,a  \\
2\,\,3\,\,4        & & 2\,\,3\,\,4\,\,5\,\,6 \\
5                  & & b \,\, c\,\, d\\
                   & & 7
\end{array}$$
where hexadecimal notation is used. (This coincides with a case of
\cite{PM}, leading to the graph ${\mathcal Q}_d$ in Corollary 9
below. Compare also with \cite{Baden,Costa}). The advantage of
this labeling is that it can be extended to all $d$-neighborhoods
of vertices of $S$ so that in each row the labels from 1 to
$q=2d^2+2d+1$ appear subsequently and contiguously from left to
right, which in the cases of Figures 6 and 7 makes that the first
row of the depicted cutout of $C_q\times C_q$ receives the
following labels from left to right:
$3\,\,4\,\,5\,\,1\,\,2\,\,3$ and 
$4\,\,5\,\,6\,\,7\,\,8\,\,9\,\,a\,\,b\,\,c\,\,d\,\,1\,\,2\,\,3\,\,4$
respectively, and so on for the remaining rows of the cutout,
comprising all the numbers in the range $1,2,\ldots, 2d(d+1)$ in
increasing order mod $q$ from left to right, from certain vertex
on, and then extended periodically for the vertices of each
horizontal path in ${\Lambda}$. Observe that the disposition of
these labels yields a $q\times q$ Latin square and that the vertices
having any fixed label constitute a translate of the
sublattice of ${\Lambda}$ generated by $\{(d,d+1),(d+1,-d)\}$.

\begin{thm}
There exists a toroidal graph $C_m\times C_n$ having a $d$-perfect code partition
${\mathcal S}_{d-1}=\{S_{d-1}^0=S_{d-1},S_{d-1}^1,\ldots,S_{d-1}^{q-1}\}$,
where $q=2d^2+2d+1$, if and only if $m$ and $n$ are multiples of $q$.
Each of the component codes $S_{d-1}^i$ of ${\mathcal S}_{d-1}$ has cardinality
$mn/q$ and cannot be obtained by side identifications from 1-perfect codes in
any rectangular grid graph.
\end{thm}

\proof
Since a cutout square $P_q\times P_q$ of
$C_q\times C_q$ projects onto
$q^2$ different vertices of $C_q\times C_q$ and since any $d$-neighborhood of a vertex of
${\Lambda}$ contains $q$ vertices, then the projected code from
$S$ in $C_q\times C_q$ contains $q$ vertices.
The statement of the theorem follows by forming larger cutouts which are rectangular
arrangements obtained by continuation of copies of the original cutout.
\qfd

\begin{cor}
Identification of vertices in ${\Lambda}$ with each common label
$i=1,2,\ldots,$ $q=2d^2+2d+1$ yields the quotient graph ${\mathcal Q}_d$
of its partition into $q$ $d$-perfect codes.
Moreover, ${\mathcal Q}_d$ is a 4-regular bipartite graph,
namely the undirected Cayley graph
of the cyclic group $Z_q$ under the generator set $\{1,2d^2\}$.
\end{cor}

\proof
What is the difference mod $q$ between the label of any vertex $w$ of ${\Lambda}$ and
the label of the vertex immediately below $w$? We answer this on the central,
longest, vertical path, or column of vertices, of $N(v)$.
Observe that the labels in this path start in its top vertex with label 1 and continue
with labels that increase
mod $q$ in subsequent increments of value $2d^2$. This is due to the fact that
by starting with the vertex labeled 1 and descending in this column
cyclically mod $2d+1$ by jumping $d$ rows in each of $2d+1$ instances,
we cover subsequently the central
position in the rows $1,d+1,2d+1,d,2d,d-1,2d-1,\ldots,2,2+d$ and again 1. By finishing this
process in label 2, two steps before returning to label 1, we
obtain the sequence
of $2d-1$ subsequent label differences $d+2,d+2,d+1,\ldots,d+2,d+1$, which starts with
$d+1$ twice and then alternating $d+1$ with $d+2$, amounting to $d$ times $d+1$ plus
$d-1$ times $d$. We conclude that
${\mathcal Q}_d$ is the 4-regular bipartite graph whose vertex
set is the cyclic group $Z_q=\{1,2,\ldots,q-1,q=0\}$, with vertices $i$ and $j$
adjacent if and only if $i-j\in\{1,-1,2d^2,-2d^2\}$. Thus, ${\mathcal Q}_d$ is the claimed
Cayley graph.
\qfd

Of course, Corollary 9 holds with ${\Lambda}$ replaced by any $C_m\times C_n$ as in
Theorem 8.
Continuing the remark that ends the second paragraph of this section, we note that
there are only
two possible partitions of ${\Lambda}$ into copies of a $d$-perfect code
of ${\Lambda}$, corresponding to the two enantiomorphic presentations of such a code.
In the second such code, the labeling of the vertices of an $N(v)$ ascends on each row
from right to left, instead of being from left to right.

\section{Total perfect code partitions of ${\Lambda}$}

\begin{figure}
\unitlength=0.30mm
\special{em:linewidth 0.4pt}
\linethickness{0.4pt}
\begin{picture}(357.39,82.62)
\put(166.39,1.62){\circle{2.00}}
\put(176.39,1.62){\circle{2.00}}
\put(206.39,1.62){\circle{2.00}}
\put(216.39,1.62){\circle{2.00}}
\put(146.39,11.62){\circle{2.00}}
\put(156.39,11.62){\circle{2.00}}
\put(166.39,11.62){\circle{2.00}}
\put(176.39,11.62){\circle{2.00}}
\put(186.39,11.62){\circle{2.00}}
\put(196.39,11.62){\circle{2.00}}
\put(206.39,11.62){\circle{2.00}}
\put(216.39,11.62){\circle{2.00}}
\put(226.39,11.62){\circle{2.00}}
\put(236.39,11.62){\circle{2.00}}
\put(146.39,21.62){\circle{2.00}}
\put(156.39,21.62){\circle{2.00}}
\put(186.39,21.62){\circle{2.00}}
\put(196.39,21.62){\circle{2.00}}
\put(226.39,21.62){\circle{2.00}}
\put(236.39,21.62){\circle{2.00}}
\put(146.39,31.62){\circle{2.00}}
\put(156.39,31.62){\circle{2.00}}
\put(166.39,31.62){\circle{2.00}}
\put(176.39,31.62){\circle{2.00}}
\put(186.39,31.62){\circle{2.00}}
\put(196.39,31.62){\circle{2.00}}
\put(206.39,31.62){\circle{2.00}}
\put(216.39,31.62){\circle{2.00}}
\put(226.39,31.62){\circle{2.00}}
\put(236.39,31.62){\circle{2.00}}
\put(166.39,41.62){\circle{2.00}}
\put(176.39,41.62){\circle{2.00}}
\put(206.39,41.62){\circle{2.00}}
\put(216.39,41.62){\circle{2.00}}
\put(146.39,51.62){\circle{2.00}}
\put(156.39,51.62){\circle{2.00}}
\put(166.39,51.62){\circle{2.00}}
\put(176.39,51.62){\circle{2.00}}
\put(186.39,51.62){\circle{2.00}}
\put(196.39,51.62){\circle{2.00}}
\put(206.39,51.62){\circle{2.00}}
\put(216.39,51.62){\circle{2.00}}
\put(226.39,51.62){\circle{2.00}}
\put(236.39,51.62){\circle{2.00}}
\put(146.39,61.62){\circle{2.00}}
\put(156.39,61.62){\circle{2.00}}
\put(186.39,61.62){\circle{2.00}}
\put(196.39,61.62){\circle{2.00}}
\put(226.39,61.62){\circle{2.00}}
\put(236.39,61.62){\circle{2.00}}
\put(146.39,71.62){\circle{2.00}}
\put(156.39,71.62){\circle{2.00}}
\put(166.39,71.62){\circle{2.00}}
\put(176.39,71.62){\circle{2.00}}
\put(186.39,71.62){\circle{2.00}}
\put(196.39,71.62){\circle{2.00}}
\put(206.39,71.62){\circle{2.00}}
\put(216.39,71.62){\circle{2.00}}
\put(226.39,71.62){\circle{2.00}}
\put(236.39,71.62){\circle{2.00}}
\put(166.39,81.62){\circle{2.00}}
\put(176.39,81.62){\circle{2.00}}
\put(206.39,81.62){\circle{2.00}}
\put(216.39,81.62){\circle{2.00}}
\put(186.39,41.62){\circle*{2.00}}
\put(196.39,41.62){\circle*{2.00}}
\put(176.39,61.62){\circle*{2.00}}
\put(166.39,61.62){\circle*{2.00}}
\put(156.39,81.62){\circle*{2.00}}
\put(146.39,81.62){\circle*{2.00}}
\put(206.39,21.62){\circle*{2.00}}
\put(216.39,21.62){\circle*{2.00}}
\put(226.39,1.62){\circle*{2.00}}
\put(236.39,1.62){\circle*{2.00}}
\put(226.39,41.62){\circle*{2.00}}
\put(236.39,41.62){\circle*{2.00}}
\put(216.39,61.62){\circle*{2.00}}
\put(206.39,61.62){\circle*{2.00}}
\put(196.39,81.62){\circle*{2.00}}
\put(186.39,81.62){\circle*{2.00}}
\put(176.39,21.62){\circle*{2.00}}
\put(166.39,21.62){\circle*{2.00}}
\put(156.39,1.62){\circle*{2.00}}
\put(146.39,1.62){\circle*{2.00}}
\put(186.39,1.62){\circle*{2.00}}
\put(196.39,1.62){\circle*{2.00}}
\put(156.39,41.62){\circle*{2.00}}
\put(146.39,41.62){\circle*{2.00}}
\emline{176.39}{23.62}{1}{176.39}{24.62}{2}
\emline{176.39}{25.62}{3}{176.39}{26.62}{4}
\emline{176.39}{27.62}{5}{176.39}{28.62}{6}
\emline{176.39}{29.62}{7}{176.39}{30.62}{8}
\emline{178.39}{21.62}{9}{179.39}{21.62}{10}
\emline{180.39}{21.62}{11}{181.39}{21.62}{12}
\emline{182.39}{21.62}{13}{183.39}{21.62}{14}
\emline{184.39}{21.62}{15}{185.39}{21.62}{16}
\put(226.39,81.62){\circle*{2.00}}
\put(236.39,81.62){\circle*{2.00}}
\emline{188.39}{21.62}{17}{189.39}{21.62}{18}
\emline{190.39}{21.62}{19}{191.39}{21.62}{20}
\emline{192.39}{21.62}{21}{193.39}{21.62}{22}
\emline{194.39}{21.62}{23}{195.39}{21.62}{24}
\emline{198.39}{21.62}{25}{199.39}{21.62}{26}
\emline{200.39}{21.62}{27}{201.39}{21.62}{28}
\emline{202.39}{21.62}{29}{203.39}{21.62}{30}
\emline{204.39}{21.62}{31}{205.39}{21.62}{32}
\emline{178.39}{61.62}{33}{179.39}{61.62}{34}
\emline{180.39}{61.62}{35}{181.39}{61.62}{36}
\emline{182.39}{61.62}{37}{183.39}{61.62}{38}
\emline{184.39}{61.62}{39}{185.39}{61.62}{40}
\emline{188.39}{61.62}{41}{189.39}{61.62}{42}
\emline{190.39}{61.62}{43}{191.39}{61.62}{44}
\emline{192.39}{61.62}{45}{193.39}{61.62}{46}
\emline{194.39}{61.62}{47}{195.39}{61.62}{48}
\emline{198.39}{61.62}{49}{199.39}{61.62}{50}
\emline{200.39}{61.62}{51}{201.39}{61.62}{52}
\emline{202.39}{61.62}{53}{203.39}{61.62}{54}
\emline{204.39}{61.62}{55}{205.39}{61.62}{56}
\emline{176.39}{33.62}{57}{176.39}{34.62}{58}
\emline{176.39}{35.62}{59}{176.39}{36.62}{60}
\emline{176.39}{37.62}{61}{176.39}{38.62}{62}
\emline{176.39}{39.62}{63}{176.39}{40.62}{64}
\emline{176.39}{43.62}{65}{176.39}{44.62}{66}
\emline{176.39}{45.62}{67}{176.39}{46.62}{68}
\emline{176.39}{47.62}{69}{176.39}{48.62}{70}
\emline{176.39}{49.62}{71}{176.39}{50.62}{72}
\emline{176.39}{53.62}{73}{176.39}{54.62}{74}
\emline{176.39}{55.62}{75}{176.39}{56.62}{76}
\emline{176.39}{57.62}{77}{176.39}{58.62}{78}
\emline{176.39}{59.62}{79}{176.39}{60.62}{80}
\emline{216.39}{23.62}{81}{216.39}{24.62}{82}
\emline{216.39}{25.62}{83}{216.39}{26.62}{84}
\emline{216.39}{27.62}{85}{216.39}{28.62}{86}
\emline{216.39}{29.62}{87}{216.39}{30.62}{88}
\emline{216.39}{33.62}{89}{216.39}{34.62}{90}
\emline{216.39}{35.62}{91}{216.39}{36.62}{92}
\emline{216.39}{37.62}{93}{216.39}{38.62}{94}
\emline{216.39}{39.62}{95}{216.39}{40.62}{96}
\emline{216.39}{43.62}{97}{216.39}{44.62}{98}
\emline{216.39}{45.62}{99}{216.39}{46.62}{100}
\emline{216.39}{47.62}{101}{216.39}{48.62}{102}
\emline{216.39}{49.62}{103}{216.39}{50.62}{104}
\emline{216.39}{53.62}{105}{216.39}{54.62}{106}
\emline{216.39}{55.62}{107}{216.39}{56.62}{108}
\emline{216.39}{57.62}{109}{216.39}{58.62}{110}
\emline{216.39}{59.62}{111}{216.39}{60.62}{112}
\put(286.39,1.62){\circle{2.00}}
\put(296.39,1.62){\circle{2.00}}
\put(326.39,1.62){\circle{2.00}}
\put(336.39,1.62){\circle{2.00}}
\put(266.39,11.62){\circle{2.00}}
\put(276.39,11.62){\circle{2.00}}
\put(286.39,11.62){\circle{2.00}}
\put(296.39,11.62){\circle{2.00}}
\put(306.39,11.62){\circle{2.00}}
\put(316.39,11.62){\circle{2.00}}
\put(326.39,11.62){\circle{2.00}}
\put(336.39,11.62){\circle{2.00}}
\put(346.39,11.62){\circle{2.00}}
\put(356.39,11.62){\circle{2.00}}
\put(266.39,21.62){\circle{2.00}}
\put(276.39,21.62){\circle{2.00}}
\put(306.39,21.62){\circle{2.00}}
\put(316.39,21.62){\circle{2.00}}
\put(346.39,21.62){\circle{2.00}}
\put(356.39,21.62){\circle{2.00}}
\put(266.39,31.62){\circle{2.00}}
\put(276.39,31.62){\circle{2.00}}
\put(286.39,31.62){\circle{2.00}}
\put(296.39,31.62){\circle{2.00}}
\put(306.39,31.62){\circle{2.00}}
\put(316.39,31.62){\circle{2.00}}
\put(326.39,31.62){\circle{2.00}}
\put(336.39,31.62){\circle{2.00}}
\put(346.39,31.62){\circle{2.00}}
\put(356.39,31.62){\circle{2.00}}
\put(286.39,41.62){\circle{2.00}}
\put(296.39,41.62){\circle{2.00}}
\put(326.39,41.62){\circle{2.00}}
\put(336.39,41.62){\circle{2.00}}
\put(266.39,51.62){\circle{2.00}}
\put(276.39,51.62){\circle{2.00}}
\put(286.39,51.62){\circle{2.00}}
\put(296.39,51.62){\circle{2.00}}
\put(306.39,51.62){\circle{2.00}}
\put(316.39,51.62){\circle{2.00}}
\put(326.39,51.62){\circle{2.00}}
\put(336.39,51.62){\circle{2.00}}
\put(346.39,51.62){\circle{2.00}}
\put(356.39,51.62){\circle{2.00}}
\put(266.39,61.62){\circle{2.00}}
\put(276.39,61.62){\circle{2.00}}
\put(306.39,61.62){\circle{2.00}}
\put(316.39,61.62){\circle{2.00}}
\put(346.39,61.62){\circle{2.00}}
\put(356.39,61.62){\circle{2.00}}
\put(266.39,71.62){\circle{2.00}}
\put(276.39,71.62){\circle{2.00}}
\put(286.39,71.62){\circle{2.00}}
\put(296.39,71.62){\circle{2.00}}
\put(306.39,71.62){\circle{2.00}}
\put(316.39,71.62){\circle{2.00}}
\put(326.39,71.62){\circle{2.00}}
\put(336.39,71.62){\circle{2.00}}
\put(346.39,71.62){\circle{2.00}}
\put(356.39,71.62){\circle{2.00}}
\put(286.39,81.62){\circle{2.00}}
\put(296.39,81.62){\circle{2.00}}
\put(326.39,81.62){\circle{2.00}}
\put(336.39,81.62){\circle{2.00}}
\put(306.39,41.62){\circle*{2.00}}
\put(316.39,41.62){\circle*{2.00}}
\put(296.39,61.62){\circle*{2.00}}
\put(286.39,61.62){\circle*{2.00}}
\put(276.39,81.62){\circle*{2.00}}
\put(266.39,81.62){\circle*{2.00}}
\put(326.39,21.62){\circle*{2.00}}
\put(336.39,21.62){\circle*{2.00}}
\put(346.39,1.62){\circle*{2.00}}
\put(356.39,1.62){\circle*{2.00}}
\put(346.39,41.62){\circle*{2.00}}
\put(356.39,41.62){\circle*{2.00}}
\put(336.39,61.62){\circle*{2.00}}
\put(326.39,61.62){\circle*{2.00}}
\put(316.39,81.62){\circle*{2.00}}
\put(306.39,81.62){\circle*{2.00}}
\put(296.39,21.62){\circle*{2.00}}
\put(286.39,21.62){\circle*{2.00}}
\put(276.39,1.62){\circle*{2.00}}
\put(266.39,1.62){\circle*{2.00}}
\put(306.39,1.62){\circle*{2.00}}
\put(316.39,1.62){\circle*{2.00}}
\put(276.39,41.62){\circle*{2.00}}
\put(266.39,41.62){\circle*{2.00}}
\emline{296.39}{23.62}{113}{296.39}{24.62}{114}
\emline{296.39}{25.62}{115}{296.39}{26.62}{116}
\emline{296.39}{27.62}{117}{296.39}{28.62}{118}
\emline{296.39}{29.62}{119}{296.39}{30.62}{120}
\emline{298.39}{21.62}{121}{299.39}{21.62}{122}
\emline{300.39}{21.62}{123}{301.39}{21.62}{124}
\emline{302.39}{21.62}{125}{303.39}{21.62}{126}
\emline{304.39}{21.62}{127}{305.39}{21.62}{128}
\put(346.39,81.62){\circle*{2.00}}
\put(356.39,81.62){\circle*{2.00}}
\emline{318.39}{21.62}{129}{319.39}{21.62}{130}
\emline{320.39}{21.62}{131}{321.39}{21.62}{132}
\emline{322.39}{21.62}{133}{323.39}{21.62}{134}
\emline{324.39}{21.62}{135}{325.39}{21.62}{136}
\emline{328.39}{21.62}{137}{329.39}{21.62}{138}
\emline{330.39}{21.62}{139}{331.39}{21.62}{140}
\emline{332.39}{21.62}{141}{333.39}{21.62}{142}
\emline{334.39}{21.62}{143}{335.39}{21.62}{144}
\emline{298.39}{61.62}{145}{299.39}{61.62}{146}
\emline{300.39}{61.62}{147}{301.39}{61.62}{148}
\emline{302.39}{61.62}{149}{303.39}{61.62}{150}
\emline{304.39}{61.62}{151}{305.39}{61.62}{152}
\emline{318.39}{61.62}{153}{319.39}{61.62}{154}
\emline{320.39}{61.62}{155}{321.39}{61.62}{156}
\emline{322.39}{61.62}{157}{323.39}{61.62}{158}
\emline{324.39}{61.62}{159}{325.39}{61.62}{160}
\emline{328.39}{61.62}{161}{329.39}{61.62}{162}
\emline{330.39}{61.62}{163}{331.39}{61.62}{164}
\emline{332.39}{61.62}{165}{333.39}{61.62}{166}
\emline{334.39}{61.62}{167}{335.39}{61.62}{168}
\emline{296.39}{53.62}{169}{296.39}{54.62}{170}
\emline{296.39}{55.62}{171}{296.39}{56.62}{172}
\emline{296.39}{57.62}{173}{296.39}{58.62}{174}
\emline{296.39}{59.62}{175}{296.39}{60.62}{176}
\emline{336.39}{23.62}{177}{336.39}{24.62}{178}
\emline{336.39}{25.62}{179}{336.39}{26.62}{180}
\emline{336.39}{27.62}{181}{336.39}{28.62}{182}
\emline{336.39}{29.62}{183}{336.39}{30.62}{184}
\emline{336.39}{53.62}{185}{336.39}{54.62}{186}
\emline{336.39}{55.62}{187}{336.39}{56.62}{188}
\emline{336.39}{57.62}{189}{336.39}{58.62}{190}
\emline{336.39}{59.62}{191}{336.39}{60.62}{192}
\emline{207.39}{21.62}{193}{215.39}{21.62}{194}
\emline{187.39}{41.62}{195}{195.39}{41.62}{196}
\emline{207.39}{61.62}{197}{215.39}{61.62}{198}
\emline{227.39}{81.62}{199}{235.39}{81.62}{200}
\emline{187.39}{81.62}{201}{195.39}{81.62}{202}
\emline{147.39}{81.62}{203}{155.39}{81.62}{204}
\emline{167.39}{61.62}{205}{175.39}{61.62}{206}
\emline{147.39}{41.62}{207}{155.39}{41.62}{208}
\emline{227.39}{41.62}{209}{235.39}{41.62}{210}
\emline{167.39}{21.62}{211}{175.39}{21.62}{212}
\emline{187.39}{1.62}{213}{195.39}{1.62}{214}
\emline{147.39}{1.62}{215}{155.39}{1.62}{216}
\emline{227.39}{1.62}{217}{235.39}{1.62}{218}
\emline{156.39}{81.62}{219}{166.39}{61.62}{220}
\emline{166.39}{61.62}{221}{156.39}{41.62}{222}
\emline{156.39}{41.62}{223}{166.39}{21.62}{224}
\emline{166.39}{21.62}{225}{156.39}{1.62}{226}
\emline{186.39}{1.62}{227}{176.39}{21.62}{228}
\emline{176.39}{21.62}{229}{186.39}{41.62}{230}
\emline{186.39}{41.62}{231}{176.39}{61.62}{232}
\emline{176.39}{61.62}{233}{186.39}{81.62}{234}
\emline{196.39}{81.62}{235}{206.39}{61.62}{236}
\emline{206.39}{61.62}{237}{196.39}{41.62}{238}
\emline{196.39}{41.62}{239}{206.39}{21.62}{240}
\emline{206.39}{21.62}{241}{196.39}{1.62}{242}
\emline{226.39}{1.62}{243}{216.39}{21.62}{244}
\emline{216.39}{21.62}{245}{226.39}{41.62}{246}
\emline{226.39}{41.62}{247}{216.39}{61.62}{248}
\emline{216.39}{61.62}{249}{226.39}{81.62}{250}
\emline{267.39}{31.62}{251}{275.39}{31.62}{252}
\emline{277.39}{31.62}{253}{285.39}{31.62}{254}
\emline{266.39}{12.62}{255}{266.39}{20.62}{256}
\emline{266.39}{22.62}{257}{266.39}{30.62}{258}
\emline{266.39}{52.62}{259}{266.39}{60.62}{260}
\emline{266.39}{62.62}{261}{266.39}{70.62}{262}
\emline{276.39}{12.62}{263}{276.39}{20.62}{264}
\emline{276.39}{22.62}{265}{276.39}{30.62}{266}
\emline{276.39}{52.62}{267}{276.39}{60.62}{268}
\emline{276.39}{62.62}{269}{276.39}{70.62}{270}
\emline{306.39}{12.62}{271}{306.39}{20.62}{272}
\emline{306.39}{22.62}{273}{306.39}{30.62}{274}
\emline{306.39}{52.62}{275}{306.39}{60.62}{276}
\emline{306.39}{62.62}{277}{306.39}{70.62}{278}
\emline{316.39}{12.62}{279}{316.39}{20.62}{280}
\emline{316.39}{22.62}{281}{316.39}{30.62}{282}
\emline{316.39}{52.62}{283}{316.39}{60.62}{284}
\emline{316.39}{62.62}{285}{316.39}{70.62}{286}
\emline{287.39}{31.62}{287}{295.39}{31.62}{288}
\emline{297.39}{31.62}{289}{305.39}{31.62}{290}
\emline{307.39}{31.62}{291}{315.39}{31.62}{292}
\emline{317.39}{31.62}{293}{325.39}{31.62}{294}
\emline{327.39}{31.62}{295}{335.39}{31.62}{296}
\emline{337.39}{31.62}{297}{345.39}{31.62}{298}
\emline{347.39}{31.62}{299}{355.39}{31.62}{300}
\emline{267.39}{11.62}{301}{275.39}{11.62}{302}
\emline{277.39}{11.62}{303}{285.39}{11.62}{304}
\emline{287.39}{11.62}{305}{295.39}{11.62}{306}
\emline{297.39}{11.62}{307}{305.39}{11.62}{308}
\emline{307.39}{11.62}{309}{315.39}{11.62}{310}
\emline{317.39}{11.62}{311}{325.39}{11.62}{312}
\emline{327.39}{11.62}{313}{335.39}{11.62}{314}
\emline{337.39}{11.62}{315}{345.39}{11.62}{316}
\emline{347.39}{11.62}{317}{355.39}{11.62}{318}
\emline{267.39}{71.62}{319}{275.39}{71.62}{320}
\emline{277.39}{71.62}{321}{285.39}{71.62}{322}
\emline{287.39}{71.62}{323}{295.39}{71.62}{324}
\emline{297.39}{71.62}{325}{305.39}{71.62}{326}
\emline{307.39}{71.62}{327}{315.39}{71.62}{328}
\emline{317.39}{71.62}{329}{325.39}{71.62}{330}
\emline{327.39}{71.62}{331}{335.39}{71.62}{332}
\emline{337.39}{71.62}{333}{345.39}{71.62}{334}
\emline{347.39}{71.62}{335}{355.39}{71.62}{336}
\emline{267.39}{51.62}{337}{275.39}{51.62}{338}
\emline{277.39}{51.62}{339}{285.39}{51.62}{340}
\emline{287.39}{51.62}{341}{295.39}{51.62}{342}
\emline{297.39}{51.62}{343}{305.39}{51.62}{344}
\emline{307.39}{51.62}{345}{315.39}{51.62}{346}
\emline{317.39}{51.62}{347}{325.39}{51.62}{348}
\emline{327.39}{51.62}{349}{335.39}{51.62}{350}
\emline{337.39}{51.62}{351}{345.39}{51.62}{352}
\emline{347.39}{51.62}{353}{355.39}{51.62}{354}
\emline{346.39}{12.62}{355}{346.39}{20.62}{356}
\emline{346.39}{22.62}{357}{346.39}{30.62}{358}
\emline{346.39}{52.62}{359}{346.39}{60.62}{360}
\emline{346.39}{62.62}{361}{346.39}{70.62}{362}
\emline{356.39}{12.62}{363}{356.39}{20.62}{364}
\emline{356.39}{22.62}{365}{356.39}{30.62}{366}
\emline{356.39}{52.62}{367}{356.39}{60.62}{368}
\emline{356.39}{62.62}{369}{356.39}{70.62}{370}
\emline{267.39}{61.62}{371}{275.39}{61.62}{372}
\emline{267.39}{21.62}{373}{275.39}{21.62}{374}
\emline{307.39}{61.62}{375}{315.39}{61.62}{376}
\emline{307.39}{21.62}{377}{315.39}{21.62}{378}
\emline{347.39}{61.62}{379}{355.39}{61.62}{380}
\emline{347.39}{21.62}{381}{355.39}{21.62}{382}
\emline{327.39}{1.62}{383}{335.39}{1.62}{384}
\emline{326.39}{2.62}{385}{326.39}{10.62}{386}
\emline{336.39}{2.62}{387}{336.39}{10.62}{388}
\emline{287.39}{1.62}{389}{295.39}{1.62}{390}
\emline{286.39}{2.62}{391}{286.39}{10.62}{392}
\emline{296.39}{2.62}{393}{296.39}{10.62}{394}
\emline{326.39}{32.62}{395}{326.39}{40.62}{396}
\emline{326.39}{42.62}{397}{326.39}{50.62}{398}
\emline{336.39}{32.62}{399}{336.39}{40.62}{400}
\emline{336.39}{42.62}{401}{336.39}{50.62}{402}
\emline{327.39}{41.62}{403}{335.39}{41.62}{404}
\emline{286.39}{32.62}{405}{286.39}{40.62}{406}
\emline{286.39}{42.62}{407}{286.39}{50.62}{408}
\emline{296.39}{32.62}{409}{296.39}{40.62}{410}
\emline{296.39}{42.62}{411}{296.39}{50.62}{412}
\emline{287.39}{41.62}{413}{295.39}{41.62}{414}
\emline{327.39}{81.62}{415}{335.39}{81.62}{416}
\emline{326.39}{72.62}{417}{326.39}{80.62}{418}
\emline{336.39}{72.62}{419}{336.39}{80.62}{420}
\emline{287.39}{81.62}{421}{295.39}{81.62}{422}
\emline{286.39}{72.62}{423}{286.39}{80.62}{424}
\emline{296.39}{72.62}{425}{296.39}{80.62}{426}
\end{picture}
\caption{Total perfect code in ${\Lambda}$ and a cutout for $C_4\times C_4$}
\end{figure}
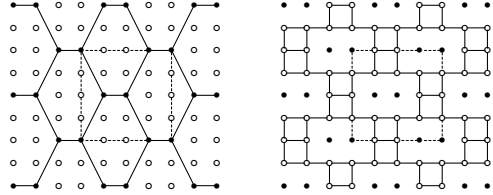

Among the PTPCs in ${\Lambda}$ arising from Theorem 1,
there is only one, call it $S_2$, such that ${\Lambda}$
admits a partition into copies of $S_2$ that can be projected
onto a partition of $C_4\times C_4$, based on
the PTPC corresponding
to the null $\{0,1\}$-sequence, and due to Corollary 5.
The total number
of copies of such PTPC $S_2$ existing in ${\Lambda}$ is 16,
of which eight have just horizontal
(vertical) induced edges, yielding four possible partitions of ${\Lambda}$ into
copies of $S_2$.

From Corollary 5,
it can be seen that there are PTPCs in toroidal graphs
$C_{4k}\times C_{4\ell}$, for $0<k,\ell\in\Z$,
obtained from $S_2$ by means of identifications in ${\Lambda}$ and
having cardinality $4k\ell$.
An example of such PTPC in $C_4\times C_4$ can be visualized in
Figure 8, where the dotted lines delineate the boundary of a cutout
from which the (flat) torus involved can be obtained; the left side of the
figure represents part of the restriction of $({\Lambda})^3$
to $S_2$; the right side represents
the complementary graph of $S_2$ in ${\Lambda}$.

\begin{thm} Let $1<m,n\in\Z$. There exists a toroidal graph $C_m\times C_n$
having a PTPC partition ${\mathcal S}_2=\{S_2^0=S_2,S_2^1,S_2^2,S_2^3\}$
if and only if $m$ and $n$ are multiples of 4.
Each component PTPC $S_2^i$ in ${\mathcal S}_2$
has cardinality $mn/4$ and cannot be obtained by side identifications
from PTPCs in any rectangular grid graph.
\qfd
\end{thm}

\section{A Penrose-tiling-like total perfect code} 

In \cite{KG}, it was shown that an $m\times n$ grid graph ${\Gamma}_{m,n}$ with
$min\{m,n\}>1$ contains a TPC if and only if
$m\equiv 0$ (mod 2) and $n\equiv -3,-1$ or 1 (mod $m+1$).

In \cite{DD2}, this was used to show that there is only one TPC $S_1$
in ${\Lambda}$ that restricts to TPCs in rectangular
grid graphs ${\Gamma}_{m,n}$, where $m$ and $n$ are integers $>2$. Moreover, the
complement ${\Lambda}\setminus S_1$ yields an aperiodic tiling of the plane (like the
Penrose tiling, \cite{Penrose}) whose automorphism group coincides with the group $D_8$ of the
square $[-\frac{1}{2},\frac{1}{2}]\times[-\frac{1}{2},\frac{1}{2}]$.

Again, this result is in contrast with the uncountability of TPCs shown in Theorem 1.

It is impossible to partition ${\Lambda}$ into copies of the graph $S_1$,
nor there are quotient toroidal graphs in ${\Lambda}$ containing a TPC obtained
by projecting $S_1$, because of the presence of a unique central ladder of area 3 in the
complementary graph of $S_1$ in ${\Lambda}$, while all the other ladders have area 2.
See Figure 9, where two concentric
stages in the construction of $S_1$ and of ${\Lambda}\setminus S_1$
are shown. Observe that the leftmost stage must be
rotated 90 degrees in order for its immersion into the second stage to be visualized,
with intermittent lines added for ease of comprehension of the immersion
(and in accordance with the definition of PDS-arrays of graphs ${\Gamma}_{m,n}$
in \cite{DD2}).

\section{PDSs in other toroidal grid graphs}

\begin{figure}
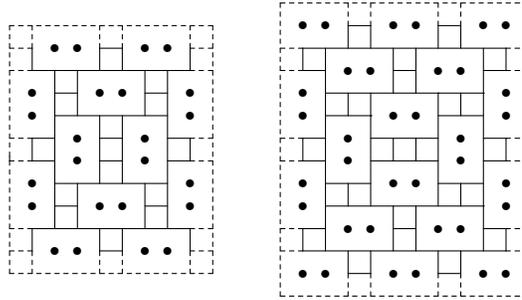

\unitlength=0.60mm
\special{em:linewidth 0.4pt}
\linethickness{0.4pt}

\caption{Two stages in the construction of $S_1$ and ${\Lambda}\setminus S_1$}
\end{figure}

Denote by $C_2$ the multigraph composed by two vertices and two parallel edges between them. We have the following complementary results.

\begin{cor}
A PTPC in a $C_2\times C_n$ exists if and only if $n$ is a multiple of 6.
\end{cor}

\proof
Let $C_2\times C_{6k}$ be represented in a (flat) torus obtained from the immersion of
the rectangular sub-grid of ${\Lambda}$ that has vertex set
$\{(i,j);0\leq i\leq 2; 0\leq j\leq 6k\}$ into the Euclidean plane,
by means of vertex identifications $(i,0)\equiv(i,6k)$,
for $0\leq i\leq 2$, and $(0,j)\equiv(2,j)$, for $0\leq j\leq 6k$, and corresponding
edge identifications. Then,
a PTPC in $C_2\times C_{6k}$ is given by the projection of the vertices
$(0,1+6i),(0,2+6i),(1,4+6i),(1,5+6i)$, where $0\leq i<k$.
\qfd

\begin{thm} There exists a toroidal graph $C_2\times C_n$
having a 1-perfect code partition ${\mathcal S}_0$
if and only if $n$ is divisible by 4.
In this case, ${\mathcal S}_0$ contains four 1-perfect codes, and
the component codes cannot be obtained by side identifications from
1-perfect codes in any rectangular grid graph.
\end{thm}

\proof Let $C_2\times C_{4k}$ be represented in a (flat) torus
obtained from the immersion of the rectangular sub-grid in
${\Lambda}$ possessing vertex set $\{(i,j);0\leq i\leq 2; 0\leq
j\leq 4k\}$ into the Euclidean plane by means of the vertex
identifications $(i,0)\equiv(i,4k)$, for $0\leq i\leq 2$, and
$(0,j)\equiv(2,j)$, for $0\leq j\leq 4k$, accompanied by the
corresponding edge identifications. Then, one of the 1-perfect codes
in ${\mathcal S}_0$ is given by the projection of the vertices
$(1,0+4i),(0,2+4i)$, where $0\leq i<k$. (An example of this for
$k=1$ is given on the left side of Figure 10, commented in a remark
below). The other three 1-perfect codes here are obtained by
translation along the vectors $(1,0),(0,1),(1,1)$. \qfd

Theorem 2.5 of \cite{LS} announces correctly the existence of the
code $S_0$ mentioned in Section 3 above. However, it also claims the
existence of a 1-perfect code in $C_4\times C_6$, which is
incorrect. The right side of Figure 10 serves to get a
counterexample: we can select successively vertices $u,v,w$ as
members of a candidate 1-perfect code in 
$C_4\times C_6$, as depicted in the figure, (where diagonals join
the dominated vertices of each of $u,v,w$; notice that after
selecting vertex $u$, the only dominating vertex for $s$, up to
symmetry, is $v$, and then the only dominating vertex for $t$ is
$w$). But then vertex $x$ cannot form part of any such 1-perfect
code, nor can be dominated by any vertex of it.

The left side of Figure 7 shows a detachment of $C_2\times C_4$
showing a 1-perfect code $\{u,v\}$, as correctly cited in the mentioned theorem,
Generally, diagonals joining dominated vertices
would form square rhombuses, but in this case,
the rhombus around $u$ is a degenerate one, because two opposite vertices of the rhombus are
identified.

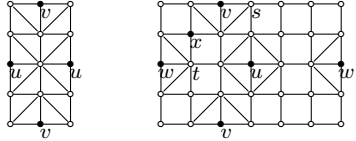
\begin{figure}
\unitlength=0.40mm
\special{em:linewidth 0.4pt}
\linethickness{0.4pt}
\begin{picture}(245.00,45.00)
\put(183.00,4.00){\circle{2.00}}
\put(193.00,4.00){\circle{2.00}}
\emline{184.00}{4.00}{1}{192.00}{4.00}{2}
\put(203.00,4.00){\circle{2.00}}
\emline{194.00}{4.00}{3}{202.00}{4.00}{4}
\put(213.00,4.00){\circle{2.00}}
\emline{204.00}{4.00}{5}{212.00}{4.00}{6}
\put(223.00,4.00){\circle{2.00}}
\emline{214.00}{4.00}{7}{222.00}{4.00}{8}
\put(233.00,4.00){\circle{2.00}}
\emline{224.00}{4.00}{9}{232.00}{4.00}{10}
\put(243.00,4.00){\circle{2.00}}
\emline{234.00}{4.00}{11}{242.00}{4.00}{12}
\put(183.00,14.00){\circle{2.00}}
\put(193.00,14.00){\circle{2.00}}
\emline{184.00}{14.00}{13}{192.00}{14.00}{14}
\put(203.00,14.00){\circle{2.00}}
\emline{194.00}{14.00}{15}{202.00}{14.00}{16}
\put(213.00,14.00){\circle{2.00}}
\emline{204.00}{14.00}{17}{212.00}{14.00}{18}
\put(223.00,14.00){\circle{2.00}}
\emline{214.00}{14.00}{19}{222.00}{14.00}{20}
\put(233.00,14.00){\circle{2.00}}
\emline{224.00}{14.00}{21}{232.00}{14.00}{22}
\put(243.00,14.00){\circle{2.00}}
\emline{234.00}{14.00}{23}{242.00}{14.00}{24}
\put(183.00,24.00){\circle{2.00}}
\put(193.00,24.00){\circle{2.00}}
\emline{184.00}{24.00}{25}{192.00}{24.00}{26}
\put(203.00,24.00){\circle{2.00}}
\emline{194.00}{24.00}{27}{202.00}{24.00}{28}
\put(213.00,24.00){\circle{2.00}}
\emline{204.00}{24.00}{29}{212.00}{24.00}{30}
\put(223.00,24.00){\circle{2.00}}
\emline{214.00}{24.00}{31}{222.00}{24.00}{32}
\put(233.00,24.00){\circle{2.00}}
\emline{224.00}{24.00}{33}{232.00}{24.00}{34}
\put(243.00,24.00){\circle{2.00}}
\emline{234.00}{24.00}{35}{242.00}{24.00}{36}
\put(183.00,34.00){\circle{2.00}}
\put(193.00,34.00){\circle{2.00}}
\emline{184.00}{34.00}{37}{192.00}{34.00}{38}
\put(203.00,34.00){\circle{2.00}}
\emline{194.00}{34.00}{39}{202.00}{34.00}{40}
\put(213.00,34.00){\circle{2.00}}
\emline{204.00}{34.00}{41}{212.00}{34.00}{42}
\put(223.00,34.00){\circle{2.00}}
\emline{214.00}{34.00}{43}{222.00}{34.00}{44}
\put(233.00,34.00){\circle{2.00}}
\emline{224.00}{34.00}{45}{232.00}{34.00}{46}
\put(243.00,34.00){\circle{2.00}}
\emline{234.00}{34.00}{47}{242.00}{34.00}{48}
\put(183.00,44.00){\circle{2.00}}
\put(193.00,44.00){\circle{2.00}}
\emline{184.00}{44.00}{49}{192.00}{44.00}{50}
\put(203.00,44.00){\circle{2.00}}
\emline{194.00}{44.00}{51}{202.00}{44.00}{52}
\put(213.00,44.00){\circle{2.00}}
\emline{204.00}{44.00}{53}{212.00}{44.00}{54}
\put(223.00,44.00){\circle{2.00}}
\emline{214.00}{44.00}{55}{222.00}{44.00}{56}
\put(233.00,44.00){\circle{2.00}}
\emline{224.00}{44.00}{57}{232.00}{44.00}{58}
\put(243.00,44.00){\circle{2.00}}
\emline{234.00}{44.00}{59}{242.00}{44.00}{60}
\emline{183.00}{5.00}{61}{183.00}{13.00}{62}
\emline{183.00}{15.00}{63}{183.00}{23.00}{64}
\emline{183.00}{25.00}{65}{183.00}{33.00}{66}
\emline{183.00}{35.00}{67}{183.00}{43.00}{68}
\emline{193.00}{5.00}{69}{193.00}{13.00}{70}
\emline{193.00}{15.00}{71}{193.00}{23.00}{72}
\emline{193.00}{25.00}{73}{193.00}{33.00}{74}
\emline{193.00}{35.00}{75}{193.00}{43.00}{76}
\emline{203.00}{5.00}{77}{203.00}{13.00}{78}
\emline{203.00}{15.00}{79}{203.00}{23.00}{80}
\emline{203.00}{25.00}{81}{203.00}{33.00}{82}
\emline{203.00}{35.00}{83}{203.00}{43.00}{84}
\emline{213.00}{5.00}{85}{213.00}{13.00}{86}
\emline{213.00}{15.00}{87}{213.00}{23.00}{88}
\emline{213.00}{25.00}{89}{213.00}{33.00}{90}
\emline{213.00}{35.00}{91}{213.00}{43.00}{92}
\emline{223.00}{5.00}{93}{223.00}{13.00}{94}
\emline{223.00}{15.00}{95}{223.00}{23.00}{96}
\emline{223.00}{25.00}{97}{223.00}{33.00}{98}
\emline{223.00}{35.00}{99}{223.00}{43.00}{100}
\emline{233.00}{5.00}{101}{233.00}{13.00}{102}
\emline{233.00}{15.00}{103}{233.00}{23.00}{104}
\emline{233.00}{25.00}{105}{233.00}{33.00}{106}
\emline{233.00}{35.00}{107}{233.00}{43.00}{108}
\emline{243.00}{5.00}{109}{243.00}{13.00}{110}
\emline{243.00}{15.00}{111}{243.00}{23.00}{112}
\emline{243.00}{25.00}{113}{243.00}{33.00}{114}
\emline{243.00}{35.00}{115}{243.00}{43.00}{116}
\put(133.00,4.00){\circle{2.00}}
\put(143.00,4.00){\circle{2.00}}
\emline{134.00}{4.00}{117}{142.00}{4.00}{118}
\put(153.00,4.00){\circle{2.00}}
\emline{144.00}{4.00}{119}{152.00}{4.00}{120}
\put(133.00,14.00){\circle{2.00}}
\put(143.00,14.00){\circle{2.00}}
\emline{134.00}{14.00}{121}{142.00}{14.00}{122}
\put(153.00,14.00){\circle{2.00}}
\emline{144.00}{14.00}{123}{152.00}{14.00}{124}
\put(133.00,24.00){\circle{2.00}}
\put(143.00,24.00){\circle{2.00}}
\emline{134.00}{24.00}{125}{142.00}{24.00}{126}
\put(153.00,24.00){\circle{2.00}}
\emline{144.00}{24.00}{127}{152.00}{24.00}{128}
\put(133.00,34.00){\circle{2.00}}
\put(143.00,34.00){\circle{2.00}}
\emline{134.00}{34.00}{129}{142.00}{34.00}{130}
\put(153.00,34.00){\circle{2.00}}
\emline{144.00}{34.00}{131}{152.00}{34.00}{132}
\put(133.00,44.00){\circle{2.00}}
\put(143.00,44.00){\circle{2.00}}
\emline{134.00}{44.00}{133}{142.00}{44.00}{134}
\put(153.00,44.00){\circle{2.00}}
\emline{144.00}{44.00}{135}{152.00}{44.00}{136}
\emline{133.00}{5.00}{137}{133.00}{13.00}{138}
\emline{133.00}{15.00}{139}{133.00}{23.00}{140}
\emline{133.00}{25.00}{141}{133.00}{33.00}{142}
\emline{133.00}{35.00}{143}{133.00}{43.00}{144}
\emline{143.00}{5.00}{145}{143.00}{13.00}{146}
\emline{143.00}{15.00}{147}{143.00}{23.00}{148}
\emline{143.00}{25.00}{149}{143.00}{33.00}{150}
\emline{143.00}{35.00}{151}{143.00}{43.00}{152}
\emline{153.00}{5.00}{153}{153.00}{13.00}{154}
\emline{153.00}{15.00}{155}{153.00}{23.00}{156}
\emline{153.00}{25.00}{157}{153.00}{33.00}{158}
\emline{153.00}{35.00}{159}{153.00}{43.00}{160}
\put(213.00,24.00){\circle*{2.00}}
\emline{214.00}{33.00}{161}{222.00}{25.00}{162}
\emline{222.00}{23.00}{163}{214.00}{15.00}{164}
\emline{212.00}{15.00}{165}{204.00}{23.00}{166}
\emline{204.00}{25.00}{167}{212.00}{33.00}{168}
\emline{212.00}{5.00}{169}{204.00}{13.00}{170}
\emline{202.00}{13.00}{171}{194.00}{5.00}{172}
\emline{212.00}{43.00}{173}{204.00}{35.00}{174}
\emline{202.00}{35.00}{175}{194.00}{43.00}{176}
\put(203.00,44.00){\circle*{2.00}}
\put(215.00,21.00){\makebox(0,0)[cc]{$_u$}}
\put(205.00,41.00){\makebox(0,0)[cc]{$_v$}}
\put(203.00,4.00){\circle*{2.00}}
\put(205.00,1.00){\makebox(0,0)[cc]{$_v$}}
\put(183.00,24.00){\circle*{2.00}}
\put(243.00,24.00){\circle*{2.00}}
\put(185.00,21.00){\makebox(0,0)[cc]{$_w$}}
\put(245.00,21.00){\makebox(0,0)[cc]{$_w$}}
\emline{184.00}{33.00}{177}{192.00}{25.00}{178}
\emline{192.00}{23.00}{179}{184.00}{15.00}{180}
\emline{242.00}{33.00}{181}{234.00}{25.00}{182}
\emline{234.00}{23.00}{183}{242.00}{15.00}{184}
\put(143.00,44.00){\circle*{2.00}}
\put(143.00,4.00){\circle*{2.00}}
\put(133.00,24.00){\circle*{2.00}}
\put(153.00,24.00){\circle*{2.00}}
\emline{134.00}{43.00}{185}{142.00}{35.00}{186}
\emline{144.00}{35.00}{187}{152.00}{43.00}{188}
\emline{134.00}{33.00}{189}{142.00}{25.00}{190}
\emline{144.00}{25.00}{191}{152.00}{33.00}{192}
\emline{142.00}{23.00}{193}{134.00}{15.00}{194}
\emline{144.00}{23.00}{195}{152.00}{15.00}{196}
\emline{134.00}{5.00}{197}{142.00}{13.00}{198}
\emline{144.00}{13.00}{199}{152.00}{5.00}{200}
\put(155.00,21.00){\makebox(0,0)[cc]{$_u$}}
\put(145.00,41.00){\makebox(0,0)[cc]{$_v$}}
\put(145.00,1.00){\makebox(0,0)[cc]{$_v$}}
\put(135.00,21.00){\makebox(0,0)[cc]{$_u$}}
\put(193.00,34.00){\circle*{2.00}}
\put(195.00,31.00){\makebox(0,0)[cc]{$_x$}}
\put(215.00,41.00){\makebox(0,0)[cc]{$_s$}}
\put(195.00,21.00){\makebox(0,0)[cc]{$_t$}}
\end{picture}
\caption{Example and counterexample for Theorem 2.5 of \cite{LS}}
\end{figure}

\section{Case of integer lattice graphs of $\R^r$, $r\ge 2$}

Let $1<r\in\Z$. In this section,
we replace ${\Lambda}={\Lambda}^2$ by the
integer lattice graph ${\Lambda}^r$ of $\R^r$.
We extend facts and results in Section 4 to
${\Lambda}^r$ as follows, noticing that this can only be done for $d=1$ if $r>1$, in the
notation of Section 4.

\begin{thm}
There exists a linear 1-perfect code $S_0$ in ${\Lambda}^r$. Moreover, $S_0$
can be taken as the sublattice of ${\Lambda}^r$ generated by the vectors
$(1,r,0,\ldots,0),$ $(2,-1,0,\dots,0),$ $(3,0,-1,0,\ldots,0),$ $\ldots,$ $(r,0,\ldots,0,-1))$.
Furthermore, $S_0$ does not restrict to a 1-perfect code
in any $r$-dimensional parallelepiped
grid graph ${\Gamma}_{m_1,m_2,\ldots,m_r}$ with at least one of
$m_1,m_2,\ldots,m_r$ larger than $4$.
\end{thm}

\proof We consider the vertices of ${\Lambda}^r$ as centers of
$r$-dimensional cubes with unit-length edges that are parallel to
the coordinate axes of $r$-space. The collection of all such cubes
forms a tessellation of $\R^r$ which can be interpreted as a Voronoi
diagram, \cite{Costa}. Each cube in such a tessellation receives a
label in the range $1,2,\ldots,2r+1$ according to the following
rule, which extends the case $d=1$ of the argument previous to
Theorem 8 and Corollary 9, above. Let $r+1$ label the null vertex of
${\Lambda}^r$; let $r+1+i\delta$ label the vertex having all
coordinates null but for the $i$-th coordinate, which equals
$\delta=\pm 1$. Extend this labeling initialization by making the
difference $\ell(w)-\ell(v)$ of the labels $\ell(v)$ and $\ell(w)$
of two contiguous vertices $v=(v_1,v_2,\ldots,v_r)$ and
$w=(w_1,w_2,\ldots,w_r)$ equal to $\ell(w)-\ell(v)=i$, where
$w_j=v_j$ if $j\neq i$ and $w_i=v_i+1$. The vertices of ${\Lambda}$
having any fixed label $1,\ldots,2r+1$ constitute a translate of the
sublattice generated by the vectors cited in the statement. The
resulting labeling of the vertices of ${\Lambda}^r$ is as claimed.
The second assertion of the statement arises from an observation in the first paragraph of
Section 4, after our citing Figure 6.\qfd

\begin{figure}
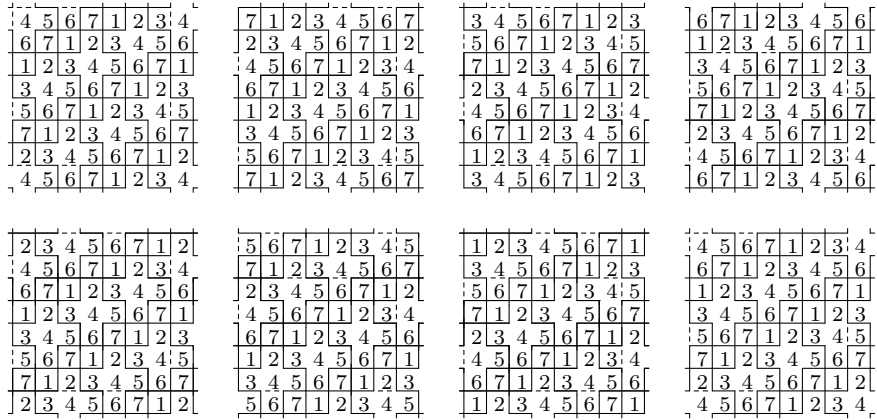

\unitlength=0.60mm
\special{em:linewidth 0.4pt}
\linethickness{0.4pt}

\caption{Case $r=3$}
\end{figure}

Figure 11 illustrates, for $r=3$, the labeling indicated in Theorem 13 by means of
a partial representation of
the mentioned Voronoi-diagram tessellation, which in this case is 3-dimensional.
The figure only shows eight partially horizontal levels of such tessellation,
where {\bf(a)} each edge represents a 4-cycle face looked upon from above; {\bf(b)}
each vertex
represents the intersection edge of two such faces; {\bf(c)} the first (upper-left) and
eighth (lower-right) levels coincide. Only faces of
outlines of 3-dimensional crosses with central cube labeled 4 and remaining six
cubes labeled 1,2,3,5,6,7 are indicated.

The centers of the $8^3$ labeled 3-cubes suggested in the figure are original vertices
of ${\Lambda}^3$.
The convex hull of these $8^3$ vertices constitutes a cutout of a 3-dimensional torus
(flat in the sense of \cite{Costa}) in which a 3-toroidal graph
$C_7\times C_7\times C_7$ is embedded, obtained by identification of (corresponding vertices
in) opposite faces of the cutout, which happens to be a 3-dimensional cube of side lengths
8, that is: each side, or maximal lateral path of the cutout is a path of length 8.
Then the distribution of the seven labels 1 through 7 yields a
1-perfect code partition of $C_7\times C_7\times C_7$, which illustrates the following
corollary. By taking the inverse image of this partition through the canonical projection
graph map $\rho:{\Lambda}^3\rightarrow C_7\times C_7\times C_7$, a 1-perfect code
partition ${\mathcal S}_0$ 
is obtained in ${\Lambda}^3$. By reverting the label ordering along one
of the coordinate directions, the enantiomorphic image of $S_0$ is obtained (through an
$(r-1)$-dimensional mirror),
as well as the enantiomorphic image of ${\mathcal S}_0$, yielding a
total of $2(2r+1)$ perfect codes equivalent to $S_0$.

Likewise, an example for ${\Lambda}^4$ leads to a 1-perfect code partition of the
4-toroidal graph $C_9\times C_9\times C_3\times C_9$, where the component $C_3$
corresponds to the third coordinate direction, with label differences
$\pm 3$ mod 9 between neighboring vertices along it.
Theorem 8 and Corollary 9 generalize now as follows.

\begin{cor} ${\Lambda}^r$ admits a partition into $2r+1$ copies of $S_0$. Moreover,
the total number of perfect codes isomorphic to $S_0$ in ${\Lambda}^r$ is $2(2r+1)$. \qfd
\end{cor}

\begin{cor}
There exists an $r$-dimensional toroidal graph
$C_{m_1}\times\ldots\times C_{m_r}$ having a 1-perfect code
partition ${\mathcal S}_0=$ $\{S_0^0=S_0,$ $S_0^1,$ $\ldots,$
$S_0^{2r}\}$ if and only if $m_i$ is a multiple of $q/\gcd(q,i)$,
for $i=1,\ldots,r$. Each component code $S_0^i$ of ${\mathcal S}_0$
has cardinality $m_1\ldots m_r/q$ and cannot be obtained by side
identifications from 1-perfect codes in any $r$-dimensional
parallelepiped grid graph. Identification of vertices in
${\Lambda}^r$ with each common label $i=1,2,\ldots, q$ yields the
quotient graph of the partition. This graph is $K_q$, that is the
undirected Cayley graph of the cyclic group $Z_q$ under the
generator set $\{1,2,\ldots,r\}$. \qfd
\end{cor}


\begin{thebibliography}{1}
\bibitem{Baden}
B. F. AlBdaiwi and M. L. Livingston,
{\it Perfect distance-d placements in 2D toroidal networks},
Jour. Supercomput., {\bf 29} (2004), 45--57.
\bibitem{BBS} D. W. Bange, A. E. Barkauskas and P. J. Slater, {\it Efficient dominating
sets in graphs}, Appl. Discrete Math, eds. R. D. Ringeisen and F. S. Roberts, SIAM,
Philadelphia, 1988, 189-199.
\bibitem{biggs} N. Biggs, Algebraic Graph Theory, Cambridge University Press, 1993.
\bibitem{Costa}S. I. Costa, M. Muniz, E. Agustini,  R. Palazzo,
{\it Graphs, tessellations, and perfect codes on flat tori},
IEEE Transact. Inform. Theory, {\bf 50} (2004), 2363--2377
\bibitem{DD2} I. J. Dejter and A. A. Delgado,
{\it Perfect domination in rectangular grid graphs},
J. Combin. Math. Combin. Comput., \textbf{70} (2009) 177--196.
\bibitem{FH} M. R. Fellows and M. N. Hoover, \textit{Perfect Domination} Australasian
J. of Combinatorics, \textbf{3} (1991), 141--150.
\bibitem{KG} W. F. Klostermeyer and J. L. Goldwasser, {\it Total Perfect Codes in
Grid Codes}, Bull. Inst. Comb. Appl., 46(2006) 61-68.
\bibitem{K} J. Kratochvil and M. Kriv\'anek, \textit{On the Computational Complexity
of Codes in Graphs}, in Proc. MFCS 1988, LN in Comp. Sci. 324
(Springer-Verlag), 396--404.
\bibitem{LS} M. Livingston and Q. F. Stout, \textit{Perfect Dominating Sets}, Congr. Numer.,
\textbf{79} (1990), 187--203.
\bibitem{Penrose}
R. Penrose, \textit{Bull. Inst. Maths. Appl.}, \textbf{10} (1974), 266.
\bibitem{W} P. M. Weichsel, \textit{Dominating Sets of n-Cubes}, J. Graph Theory,
\textbf{18} (1994), 479--488.
\bibitem{PM} J. L. A. Yebra, M. A. Fiol, P. Morillo and I. Alegre,
{\it The diameter of undirected
graphs associated to plane tessellations}, Ars Combinatoria, \textbf{20B}(1985), 159--171.
\end{thebibliography}
\end{document}